\input amstex
\documentstyle{amsppt}
\magnification 1200
\vcorrection{-1cm}
\NoBlackBoxes
\input epsf

\rightheadtext{ Counting lattice triangulations }
\topmatter
\title      
            Counting lattice triangulations:
            Fredholm equations in combinatorics
\endtitle
\author     S.~Yu.~Orevkov
\endauthor
\abstract
            Let $f(m,n)$ be the number of primitive lattice triangulations
            of $m\times n$ rectangle. We compute the limits
            $\lim_n f(m,n)^{1/n}$ for $m=2$ and $3$. For $m=2$ we obtain
            the exact value of the limit which is equal to $(611+\sqrt{73})/36$.
            For $m=3$, we express the limit in terms of certain Fredholm's
            integral equation on generating functions. This provides a
            polynomial time algorithm for computation of the limit with
            any given precision (polynomial with respect to the number
            of computed digits).
\endabstract
\address
Steklov Mathematical Institute, Gubkina 8, Moscow, Russia
\endaddress

\address
IMT, l'universit\'e Paul Sabatier, 118 route de Narbonne, Toulouse, France
\endaddress

\address
AGHA Laboratory, Moscow Institute of Physics and Technology, Russia
\endaddress

\email
orevkov\@math.ups-tlse.fr
\endemail

\endtopmatter

\def\Z{\Bbb Z}

\def\R{\Bbb R}
\def\C{\Bbb C}
\def\RP{\Bbb{RP}}
\def\Re{\operatorname{Re}}
\def\Im{\operatorname{Im}}

\def\np{{\operatorname{np}}}

\def\refAnclin  {1}
\def\refFred    {2}
\def\refGKZ     {3}
\def\refKZ      {4}
\def\refKK      {5}
\def\refKhved   {6}
\def\refMVW     {7}
\def\refOrevkov {8}
\def\refOrKh    {9}
\def\refShSh    {10} 
\def\refTamar   {11}
\def\refWelzl   {12}
\def\refWelzlTalk {13}

\def\figShad     {1}
\def\figTetrisA  {2}
\def\figT        {3}
\def\figTT       {4}
\def\figFGH      {5}
\def\figPsiLevel {6}
\def\figProg     {7}
\def\figPsiImage {8}
\def\figPsiImgII {9}

\def\tabCapacity {1}
\def\tabFive  {2}
\def\tabSix   {3}
\def\tabSeven {4}
\def\tabEight {5}
\def\tabNine  {6}
\def\tabErr {7}

\def\eqCone {1}
\def\eqF    {2}
\def\eqFF   {3}
\def\eqGG   {4}
\def\eqHH   {5}
\def\eqFFq  {6}
\def\eqGGq  {7}
\def\eqHHq  {8}
\def\eqPW   {9}
\def\eqDefP {10}
\def\eqg    {11}
\def\eqPsi  {12}
\def\eqHg   {13}

\def\eqHs   {14}
\def\eqPhi  {15}
\def\eqMinPsi   {16}
\def\eqMinP     {17}
\def\eqFredG    {18}
\def\eqFredTau  {19}
\def\eqIntH     {20}

\def\eqErrL     {21}
\def\eqFredGen  {22}
\def\eqFredHolo {23}
\def\eqErrDphi  {24}

\def\eqFredAppr   {25}
\def\eqFredErrA   {26}
\def\eqFredErrAA  {27}
\def\eqFredErrB   {28}

\def\eqFredErrRho {30}
\def\eqFredParam  {31}
\def\eqFredParamPf{32}

\def\eqNPa       {33}
\def\eqNPb       {34}
\def\eqNPc       {35}

\def\lemTetrisA {2.1}
\def\lemTetrisB {2.2}
\def\exTetrisA  {2.3}
\def\exTetrisB  {2.4}
\def\conjConv   {2.5}
\def\propConv   {2.6}

\def\lemFone  {3.1}
\def\lemFtwo  {3.2}

\def\lemP     {4.1}
\def\lemPhi   {4.2}
\def\lemG     {4.3}
\def\lemFredG {4.4}
\def\propFred {4.5}

\def\lemErrL  {5.1}
\def\lemTaylor {5.2}
\def\lemFredHolo {5.3}
\def\lemFred     {5.4}
\def\lemFredEV   {5.5}
\def\lemParam    {5.6}

\def\propNP     {6.1}

\def\sectRecursion {2}
\def\sectCtwo      {3}
\def\sectCthree    {4}
\def\sectFred      {5}
\def\sectNonprim   {6}

\document

\head 1. Introduction
\endhead
A {\it lattice triangulation} of a (lattice) polygon in $\R^2$ is a triangulation
with all vertices in $\Z^2$.
As it was discovered in [\refGKZ], lattice triangulations are important
in algebraic geometry (see also [\refOrKh]).
A lattice triangulation is called {\it primitive} (or {\it unimodular})
if each triangle is primitive, i.e., has the minimal possible area $1/2$.
We denote the number of primitive lattice triangulations of the
rectangle $m\times n$ by $f(m,n)$. Let
$$
\split
         &c(m,n) = { \log_2f(m,n)\over mn},
\qquad
         c_m = \sup_n c(m,n) = \lim_{n\to\infty} c(m,n),
\\
         &c = \sup_m c_m = \lim_{m\to\infty} c_m
            = \sup_n c(n,n) = \lim_{n\to\infty} c(n,n).
\endsplit
$$
The existence of the limits is proven in [\refKZ; Proposition 3.6]. The number
$c(m,n)$ is called in [\refKZ] the {\it capacity} of the rectangle $m\times n$.
In [\refOrevkov] I gave an upper bound $c<6$ (which can be easily improved by the same arguments up to
$c<\log_2 27= 4.755$: it is enough just not to distinguish the cases
$v_j=1$ and $v_j=2$ in the notation of [\refOrevkov]).
Later on, a much better estimate $c<3$ was obtained by Anclin [\refAnclin]
as well as $c_m<3-1/m$.
A yet better upper bound $c<4\log_2\frac{1+\sqrt5}2 = \log_2 6.854 = 2.777$ is obtained in [\refMVW] and announced in [\refWelzl]
(I have not seen the manuscript [\refMVW] but Professor Welzl kindly sent me the slides
of his talk [\refWelzlTalk] where the proof of this bound is clearly exposed).

Easy to see that
$$
   f(1,n)=\binom{2n}n  \qquad\text{whence}\quad c_1 = 2                        \eqno(\eqCone)
$$
which yields a lower bound $c>2$. It is also computed in [\refKZ] that
$c\ge c_4\ge c(4,32)=2.055702$.
It is written in [\refKZ; \S2.1]:
{\sl``For $f(2,n)$ we have no explicit formula, and we cannot evaluate the asymptotics precisely''}.
We still have no explicit formula for $f(2,n)$ but we give here the principal
term of the asymptotics:

\proclaim{ Theorem 1 } $\lim_{n\to\infty}f(2,n)^{1/n}=\alpha$ where
$$
   \alpha={611+\sqrt{73}\over 36},
   \quad\text{hence}\quad
    c_2 = {1\over 2}\log_2\alpha = 2.05256897
$$
\endproclaim

An exact value of $c_3$,
in a sense, is given in Proposition~\propFred\
where we express $c_3$ in terms of Fredholm's integral equations on
certain generating functions. In particular,
Proposition~\propFred\ provides an algorithm to compute $c_3$ up to $n$
digits in a polynomial time in $n$. A Mathematica code implementing
the main step of this algorithm is presented in Figure~\figProg\ below.

\proclaim{ Theorem 2 } $\lim_{n\to\infty} f(3,n)^{\frac1{3n}}$, up to $360$ digits, is equal to
$$
\split
 4.&239369481548025671877625742045235772100695711251795499830801
 \\&687833358238276728987837054831763341276708855553395893005289
 \\&580195934799338289257489707990192054275721787374165246347114
 \\&466096241741151814326914780021501337938335813142441896953051
 \\&597942032082556780952912032761797534112146994900056374798271
 \\&988378451540168358202181556482461979420039542105330977266751
\endsplit
$$
and hence $c_3=2.0838497...$
\endproclaim

We computed $c_3$ with this high precision hoping to find an
algebraic equation for it, or to relate it with some known constants,
but we did not succeed so far.

In \S\sectRecursion.2 we present the results of computations
of exact values of the numbers $f(m,n)$ for some small $m$ and $n$. 
These computations show in particular that
$c\ge c(5,115)=  2.10449551...$

In \S\sectNonprim\ we give an asymptotic upper bound for the number
of all (not necessarily primitive) lattice triangulations.
However it seems to be far from optimal.


\head\sectRecursion. Recurrent relations for strips of fixed width
\endhead

\subhead\sectRecursion.1. Recurrent relations
\endsubhead
Given a polygon $P\subset\R^2$, the {\it upper part} of its
boundary is the set $\{(x,y)\in P\mid y'>y\Rightarrow (x,y')\not\in P\}$.
A {\it vertical side} of $P$ is a side of $P$ contained in a line $\{x=x_0\}$.
Let $\Cal T$ be a triangulation of a polygon $P$ in $\R^2$.
We say that $Q$ is a {\it tile of} $\Cal T$ in the following three cases:
\roster
\item
   $Q$ is a triangle of $\Cal T$ without vertical sides;
\item
   $Q$ is a triangle of $\Cal T$ whose vertical side lies on the
   boundary of $P$; 
\item
  $Q$ is a union of two triangles of $\Cal T$ which share a common vertical side.
\endroster

A polygon is called $y$-{\it convex} if its intersection
with any line $x=\text{const}$ is either the empty set, or a point, or a segment.

\proclaim{ Lemma \lemTetrisA }
Let $\Cal T$ be a triangulation of a $y$-convex polygon $P$ in $\R^2$.
Then there exists a tile $Q$ of $\Cal T$ such that the upper part of
the boundary of $Q$ is contained in the upper part of the
boundary of $P$.
\endproclaim

\demo{ Proof }
Let $\Gamma_P$ be the upper part of the boundary of $P$.
Let $Q_1,\dots,Q_n$ be all the tiles of $\Cal T$ which have at least one side lying on 
$\Gamma_P$. Let $\Gamma_i$ be the union of the sides of $Q_i$ lying on $\Gamma_P$.
It is clear that each $\Gamma_i$ is either a side of $Q_i$ or a union of
two sides with a common vertex. It is also clear that the projections of the $\Gamma_i$ onto the $x$-axis have pairwise disjoint interiors, hence we may assume that
$\Gamma_1,\dots,\Gamma_n$ are numbered from the left to the right.
We say that a tile $Q_i$ is {\it shadowed on the left} (resp. {\it shadowed
on the right}) if the upper part of the boundary of $Q_i$ contains a segment $I$ such
that $I\not\subset\Gamma_P$ and $I$ is on the left (resp. on the right) of $\Gamma_i$;
see Figure~\figShad.
It is clear that none of the tiles $Q_1,\dots,Q_n$ can be shadowed
on the left and on the right simultaneously.
Hence, without lost of generality we may assume that at least one of these tiles is not
shadowed on the right. Let $i_0$ be the minimal number such that
$Q_{i_0}$ is not shadowed on the right. Then $Q_{i_0}$
is the desired tile with the upper part contained in $\Gamma_P$. Indeed, it is not shadowed on the right by its definition. It cannot be shadowed on the left neither
because otherwise $Q_{i_0-1}$ would not be shadowed on the right which
contradicts the minimality of $i_0$.
\qed
\enddemo

\midinsert
\centerline{\epsfxsize=50mm\epsfbox{shad.eps}}
\botcaption{Figure \figShad}
\endcaption
\endinsert

Now we fix an integer $m>0$ and we consider primitive lattice triangulations of
polygons contained in the vertical strip $\{0\le x\le m\}$
bounded by two graphs of continuous piecewise linear functions.

By analogy with the terminology introduced in [\refKZ, \S2.2],
we say that $\varphi:[0,m]\to\R$ is an {\it admissible function} if it is a
continuous piecewise linear function whose graph is a union of segments with
endpoints at $\Z^2$.
Let us fix an admissible function $\varphi_0$ and say that
a function $\varphi:[0,m]\to\R$ is  $\varphi_0$-{\it admissible} if
it is admissible and $\varphi(x)\ge \varphi_0(x)$ for any $x\in[0,m]$.
A $\varphi_0$-{\it admissible shape} is a polygon $S$ of the form
$\{(x,y)\in\R^2\mid 0\le x\le m,\; \varphi_0(x)\le y\le\varphi(x)\}$ for some
$\varphi_0$-admissible function $\varphi$.

As in the above definition of a tile of a triangulation, we say
that $Q$ is a {\it primitive lattice tile} in the following three cases:
\roster
\item
 $Q$ is a primitive lattice triangle without vertical sides;
\item
$Q$  is a primitive lattice triangle whose vertical side
is contained in the boundary of the strip $0\le x \le m$;
\item
$Q=\Delta_1\cup\Delta_2$ where $\Delta_1$ and $\Delta_2$ are primitive lattice
triangles such that $\Delta_1\cap\Delta_2$ is a common vertical side of $\Delta_1$
and $\Delta_2$.
\endroster
A primitive lattice tile $Q$ is $P$-{\it maximal} for a polygon $P$ if $Q\subset P$ and
the upper part of the boundary of $Q$ is contained in the upper part of the
boundary of $P$.
We say that $S'$ is a $\varphi_0$-{\it admissible subshape} of
a $\varphi_0$-admissible shape $S$, if $S'$ is the closure of
$S\setminus(Q_1\cup\dots\cup Q_n)$, where $Q_1,\dots,Q_n$ are
$S$-maximal primitive lattice tiles with pairwise disjoint interiors.
Following [\refKZ], in this case we set $\#(S',S)=n$.

Let us denote the number of primitive lattice triangulations of a polygon $P$ by $f^*(P)$.
When $P$ sits in the strip $\{0\le x\le m\}$,
we also define $f(P)$ as the
number of primitive lattice triangulations of $P$ which do not have
any interior edge whose projection onto the $x$-axis is the whole segment $[0,m]$
(we choose a simpler notation for a more complicated notion because the numbers $f(P)$
will be used more often than $f^*(P)$). 

The following lemma is the inclusion-exclusion formula in our setting.
The proof is the same as for [\refKZ, Lemma~2.2].

\proclaim{ Lemma \lemTetrisB } For any $\varphi_0$-admissible shape $S$,
we have
$$
    f^*(S)=\sum_{S'} (-1)^{\#(S',S)-1}f^*(S'),
\qquad\text{and}\qquad
    f(S)=\sum_{S'} (-1)^{\#(S',S)-1}f(S'),
$$
where the left sum is taken over all proper $\varphi_0$-admissible subshapes of $S$,
and the right sum is taken over those proper $\varphi_0$-admissible subshapes of $S$
whose upper part of the boundary contains a point from $\Z^2\cap\{0<x<m\}$.
\endproclaim

\medskip\noindent
{\bf Example \exTetrisA.}
Let $m=2$ and $\varphi_0=0$. For non-negative integers $a,b,c$,
let $S_{a,b,c}$ be the $\varphi_0$-admissible shape bounded from above
by the segment $[(0,a),(1,b)]$ and $[(1,b),(2,c)]$.
Let $f_{a,b,c}=f(S_{a,b,c})$. We set also $f_{a,b,c}=0$ when $\min(a,b,c)<0$.
Then (see Figure~\figTetrisA) the recurrent formula of Lemma~\lemTetrisB\ reads
$$
   f_{a,b,c} = \cases
      f_{a-1,b,c} + f_{a,b-1,c}+f_{a,b,c-1} - f_{a-1,b,c-1}
     &\text{if $(a,b,c)\ne(0,0,0)$,}\\
        1 &\text{if $(a,b,c)=(0,0,0)$.}
\endcases
$$
Let $F(x,y,z)=\sum_{a,b,c}f_{a,b,c}\,x^ay^bz^c$ be the generating function.
Then, by summating the recurrent relation over all triples $(a,b,c)\ne(0,0,0)$,
we obtain
$$
\split
F(x,y,z)-1 &= \sum f_{a-1,b,c}\,x^ay^bz^c+\sum f_{a,b-1,c}\,x^ay^bz^c+\dots
 \\&= \sum f_{a,b,c}\,x^{a+1}y^bz^c+\sum f_{a,b,c}\,x^ay^{b+1}z^c+\dots
 \\&= F(x,y,z)(x+y+z-xz)
\endsplit
$$
whence $F(x,y,z)=1/(1-x-y-z+xz)$.

\midinsert
\centerline{\epsfxsize=80mm\epsfbox{tetris1.eps}}
\botcaption{Figure \figTetrisA}
\endcaption
\endinsert

\medskip\noindent
{\bf Example \exTetrisB.}
Let $\varphi_0$ and $S_{a,b,c}$ be as in Example~\exTetrisA.
For non-negative $a,c$ such that $a\equiv c+1\mod 2$, we define
$S'(a,c)$ as the $\varphi_0$-admissible shape bounded from above by
the segment $[(0,a),(2,c)]$.
Let $f^*_{a,b,c}=f^*(S_{a,b,c})$ and $g^*(a,c)=f^*(S'_{a,c})$.
We set also $f^*_{a,b,c}=0$ when $\min(a,b,c)<0$ and $g^*(a,c)=0$ when $\min(a,c)<0$
or $a\equiv c \mod2$.
Then, for $(a,b,c)\ne(0,0,0)$, the recurrent formula of Lemma~\lemTetrisB\
applied to $S_{a,b,c}$ reads
$$
   f^*_{a,b,c} =
      f^*_{a-1,b,c} + f^*_{a,b-1,c}+f^*_{a,b,c-1} - f^*_{a-1,b,c-1}
      +\chi_{a,b,c}\,g^*_{a,c}
$$
where $\chi_{a,b,c}=1$ if $2b+1=a+c$, and $\chi_{a,b,c}=0$ otherwise.
Let $F^*(x,y,z)$ and $G^*(x,z)$ be the respective generating functions.
Then (cf.~Example~\exTetrisA) we have
$$
   F^*(x,y,z)-1 = F^*(x,y,z)(x+y+z-xz) + \sum\chi_{a,b,c}\,g^*_{a,c}\,x^ay^bz^c 
$$
and the last sum is equal to
$$
  \sum_{a,c} g^*_{a,c} x^a y^{(a+c-1)/2} z^c
   =y^{-1/2}\sum_{a,c} g^*_{a,c} (xy^{1/2})^a (y^{1/2}z)^c
   =y^{-1/2}G^*(xy^{1/2},y^{1/2}z)
$$
which gives us the relation
$$
     F^*(x,y,z)(1-x-y-z+xz) = 1 + y^{-1/2} G^*(xy^{1/2},y^{1/2}z).
$$
Now let us apply the recurrent relation to $S'_{a,c}$. The only
admissible subshape of $S'(a,c)$ is $S(a,(a+c-1)/2,c)$, hence
the relation for $S'_{a,c}$ reads $g^*_{a,c} = f^*_{a,(a+c-1)/2,c}$.
In terms of the generating functions this means that
$$
\split
  G^*(x,z) &= \sum_{a,c} f^*_{a,(a+c-1)/2,c}x^a z^c
 = {\text{coef}}_{u^0}\Big[\sum_{a,b,c} f^*_{a,b,c}x^a u^{2b-a-c+1} z^c\Big]
\\& = {\text{coef}}_{u^0}\Big[u\sum_{a,b,c} f^*_{a,b,c}(x/u)^a (u^2)^b (z/u)^c\Big]
 = {\text{coef}}_{u^0}\big(u F^*(x/u,u^2,z/u)\big).
\endsplit
$$

\subhead\sectRecursion.2. Some exact values of $f(m,n)$
\endsubhead

The recurrent relations in Lemma~\lemTetrisB\ provide an algorithm
of computation of exact values of $f(m,n)$ for small $m$ and $n$.
The algorithm is similar to the one described in [\refKZ, \S2.2].
We performed computations using this algorithm and one can see in
Table~\tabCapacity\ that we advanced much further with respect to
the computations in [\refKZ]. There are three reasons for this which
have more or less equal impact.

\midinsert
\hbox to 100mm{\hfill Table \tabCapacity}
\medskip
\centerline{
\vbox{\offinterlineskip
\def\h {height2pt&\omit&&\omit&\cr}
 \def\t{\times 10} \def\s{\;\;\;} \def\ss{\s\s} 
\def\q{\quad\;}
\hrule
\halign{&\vrule#&\strut\;\hfil#\hfil\,\cr
\h
& Capacities computed in [\refKZ] && Capacities computed in this paper &\cr\h
\h
 \noalign{\hrule}
\h
&$\;\q c_1=2.0000\quad    c_{4,32}=2.0557\;$&
                        &$\;c_1=2.0000\s c_{4,200}=2.0946\s c_{7,20}=2.0813\;$&\cr\h
&$\;c_{2,375}=2.0441\s c_{5,12}=2.0175\;$&
                        &$\;c_2=2.0526\s c_{5,115}=2.1045\s c_{8,13}=2.0669\;$&\cr\h
&$\;\;c_{3,60\;} =2.0275\s\; c_{6,7\;\,} =1.9841\;$&
                    &$\;c_3=2.0838\s c_{6,50\;} =2.1024\s\; c_{9,9\;} =2.0490\;$&\cr\h
\noalign{\hrule}
}}
}
\endinsert

The first reason (an evident one) is that the computers became more
powerful. The second reason is that we used another definition of
admissible shapes which allowed us to divide the amount of used memory
by $3^{m-1}$ which is rather important when $m=9$ (as it is pointed out in [\refKZ],
for this kind of algorithms,
{\sl``the bottleneck in the computations is always memory''}).
The third reason is that instead of long arithmetics, we used computations
mod different primes and then recovered the results with the Chinese Remainder
Theorem. This trick allowed us to ``convert'' memory to time whose lack was
not so crucial.

We have computed $f(3,n)$ till $n=600$ and $f(4,n)$ till $n=200$.
The exact value of $f(3,600)$ has 1127 digits and it yields $c_{3,600}=2.07966...$
Comparing this with the limit value $c_3=2.08385$ we see that the convergence
is very slow. For $m=4$, the last computed exact value is
$$
\split
  f(4,200)=\;
  &262199334303965073140522141167072596609151907003573304927487
\\&419128543906730659218480439253346584137204205604500628092962
\\&697997426095545403404830271634194339979807927812812142668569
\\&097560203843935394728621308903256950859658838687531965864231
\\&570521446370439565640979852878302993978768696718322811686043
\\&307749541067654061321020767838164602474781629699981105797912
\\&385346265396601164596410043968216134349971638142523003353406
\\&530183843913302635663917084864069175263416748948835535483336
\\&4717309018125451550646500; \qquad c(4,200)=2.09455...
\endsplit
$$
In Tables~\tabFive--\tabNine\ we present some other results of computations
in the same format as in [\refKZ]. All the computed exact values are available
on the webpage

\noindent
https://www.math.univ-toulouse.fr/\~{}orevkov/tr.html

\def\TL#1#2#3{\noindent\hskip0pt%
  \hbox to 15pt{\hfill#1}%
  \hbox to 305pt{\hfill#2}\hskip5pt
  \hbox to 32pt{\hfill#3}}

\midinsert
\hbox to 100mm{\hfill Table \tabFive}
\medskip
\hrule
\vskip 3pt
\TL{$n$}{$\#$ primitive triangulations of rectangle $5\times n$}{$c(5,n)$}
\vskip 3pt
\hrule
\vskip 3pt
\TL{1}{252}{1.5954}\par
\TL{2}{182132}{1.7474}\par
\TL{3}{182881520}{1.8297}\par
\TL{4}{208902766788}{1.8802}\par
\TL{5}{260420548144996}{1.9155}\par
\TL{6}{341816489625522032}{1.9415}\par
\TL{7}{464476385680935656240}{1.9615}\par
\TL{8}{645855159466371391947660}{1.9773}\par
\TL{9}{913036902513499041820702784}{1.9902}\par
\TL{10}{1306520849733616781789190513820}{2.0008}\par
\TL{11}{1887591165891651253904039432371172}{2.0098}\par
\TL{12}{2747848427721241461905176361078147168}{2.0174}\par
\TL{13}{4024758386310801427793602374466243714608}{2.0240}\par
\TL{14}{5924744736041718687622958191829471010847132}{2.0298}\par
\TL{15}{8757956199571261116690226598764501142088496860}{2.0348}\par
\TL{16}{12991215957916577635251095613859465176216530106080}{2.0394}\par
\TL{17}{19327902156972014645215931908930612218954616366464668}{2.0434}\par
\TL{18}{28828843648796117963238681180919362090157971920576213992}{2.0470}\par
\TL{$\vdots\;$}{$\vdots\;$}{$\vdots\;$}\par
\TL{115}{18700706608364882730712710491937598381242505216572196}{}\par
\TL{   }{74626658766824095096227084981348969054292582022965697}{}\par
\TL{   }{97536209347455134357618461876316197344892595460029612}{}\par
\TL{   }{59669310339853198410108464789290118181041289819323068}{}\par
\TL{   }{31435995596306245022821112218622320544399050742600358}{}\par
\TL{   }{31426475886050757674088153732325783413307209633451618}{}\par
\TL{   }{73035677107305109076541667755690839416820326596\hskip30pt}{2.1044}\par
\vskip 3pt
\hrule
\endinsert

\def\TL#1#2#3{\noindent\hskip0pt%
  \hbox to 12pt{\hfill#1}%
  \hbox to 308pt{\hfill#2}\hskip5pt
  \hbox to 32pt{\hfill#3}}

\midinsert
\hbox to 100mm{\hfill Table \tabSix}
\medskip
\hrule
\vskip 3pt
\TL{$n$}{$\#$ primitive triangulations of rectangle $6\times n$}{$c(6,n)$}
\vskip 3pt
\hrule
\vskip 3pt
\TL{1}{924}{1.6419}\par
\TL{2}{2801708}{1.7848}\par
\TL{3}{12244184472}{1.8617}\par
\TL{4}{61756221742966}{1.9088}\par
\TL{5}{341816489625522032}{1.9415}\par
\TL{6}{1999206934751133055518}{1.9655}\par
\TL{7}{12169409954141988707186052}{1.9840}\par
\TL{8}{76083336332947513655554918994}{1.9987}\par
\TL{9}{484772512167266688498399632918196}{2.0107}\par
\TL{10}{3131521959869770128138491287826065904}{2.0206}\par
\TL{11}{20443767611927599823217291769468449488548}{2.0289}\par
\TL{12}{134558550368400096364589064704536849131736024}{2.0360}\par
\TL{13}{891513898740246853038326950483812868791208442016}{2.0421}\par
\TL{14}{5938780824869668513059568892370775952933721743377354}{2.0474}\par
\TL{15}{39738456660509411434285642370153959115525603844258515860}{2.0521}\par
\TL{$\vdots\;$}{$\vdots\;$}{$\vdots\;$}\par
\TL{50}{733088849377871573475229677373109896289395791929}{}\par
\TL{  }{288892292779893207423013116473882328714681504398}{}\par
\TL{  }{803902969400882970235141773360945092837017232937}{}\par
\TL{  }{1864995986534063127990363531908201551410584718\hskip10pt}{2.1023}\par
\vskip 3pt
\hrule
\endinsert

\midinsert
\hbox to 100mm{\hfill Table \tabSeven}
\medskip
\hrule
\vskip 3pt
\TL{$n$}{$\#$ primitive triangulations  of rectangle $7\times n$}{$c(7,n)$}
\vskip 3pt
\hrule
\vskip 3pt
\TL{1}{3432}{1.6778}\par
\TL{2}{43936824}{1.8134}\par
\TL{3}{839660660268}{1.8862}\par
\TL{4}{18792896208387012}{1.9307}\par
\TL{5}{464476385680935656240}{1.9615}\par
\TL{6}{12169409954141988707186052}{1.9840}\par
\TL{7}{332633840844113103751597995920}{2.0014}\par
\TL{8}{9369363517501208819530429967280708}{2.0152}\par
\TL{9}{269621109753732518252493257828413137272}{2.0264}\par
\TL{10}{7880009979020501614060394747170100093057300}{2.0357}\par
\TL{11}{233031642883906149386619647304562977586311372556}{2.0435}\par
\TL{12}{6953609830304518024125545674642770582274167760568260}{2.0501}\par
\TL{13}{208980994833103266855771653608680330159883854051275967612}{2.0559}\par
\TL{$\vdots\;$}{$\vdots\;$}{$\vdots\;$}\par
\TL{20}{52066212145180734892042606757684021681422119}{}\par
\TL{  }{85233630730198914071476153736678384063983252}{2.0813}\par
\vskip 3pt
\hrule
\endinsert

\midinsert
\hbox to 100mm{\hfill Table \tabEight}
\medskip
\hrule
\vskip 3pt
\TL{$n$}{$\#$ primitive triangulations of rectangle  $8\times n$}{$c(8,b)$}
\vskip 3pt
\hrule
\vskip 3pt
\TL{1}{12870}{1.7064}\par
\TL{2}{698607816}{1.8362}\par
\TL{3}{58591381296256}{1.9056}\par
\TL{4}{5831528022482629710}{1.9480}\par
\TL{5}{645855159466371391947660}{1.9773}\par
\TL{6}{76083336332947513655554918994}{1.9987}\par
\TL{7}{9369363517501208819530429967280708}{2.0152}\par
\TL{8}{1191064812882685539785713745400934044308}{2.0282}\par
\TL{9}{155023302820254133629368881178138076738462112}{2.0388}\par
\TL{10}{20527337238769032315796332007167102984745417344046}{2.0476}\par
\TL{11}{2753810232976351788081274786378733309236298426977203848}{2.0550}\par
\TL{12}{373119178357778061717948099980013460229206030805799398500854}{2.0613}\par
\vskip2pt
\TL{13}{509513267535377736964009580351904}{}\par
\TL{  }{45392087069512323700346738258636\hskip5pt}{2.0668}\par
\vskip 3pt
\hrule
\endinsert

\midinsert
\hbox to 100mm{\hfill Table \tabNine}
\medskip
\hrule
\vskip 3pt
\TL{$n$}{$\#$ primitive triangulations of rectangle  $9\times n$}{$c(9,n)$}
\vskip 3pt
\hrule
\vskip 3pt
\TL{1}{48620}{1.7299}\par
\TL{2}{11224598424}{1.8547}\par
\TL{3}{4140106747178292}{1.9214}\par
\TL{4}{1835933384812941453312}{1.9621}\par
\TL{5}{913036902513499041820702784}{1.9902}\par
\TL{6}{484772512167266688498399632918196}{2.0107}\par
\TL{7}{269621109753732518252493257828413137272}{2.0264}\par
\TL{8}{155023302820254133629368881178138076738462112}{2.0388}\par
\TL{9}{91376512409462235694151119897052344522006298310908}{2.0489}\par
\vskip 3pt
\hrule
\endinsert

\subhead\sectRecursion.3. Convexity conjecture for the numbers $f(m,n)$
\endsubhead
The following conjecture is confirmed by all the computed exact values of
the numbers $f(m,n)$ (we set by convention $f(m,0)=1$).

\medskip\noindent
{\bf Conjecture \conjConv.}
One has $f(m,n-1)f(m,n+1)\ge f(m,n)^2$ for any $m,n\ge 1$.

\proclaim{ Proposition \propConv }
If Conjecture~\conjConv\ holds true, then
$c_m\ge (n+1)c(m,n+1)-nc(m,n)$ for any $m,n\ge 1$.
In particular, Conjecture~\conjConv\ would imply that
$c\ge c_{115}\ge 5c(115,5)-4c(115,4) = 2.1684837\dots$
\endproclaim

\demo{ Proof } Let us set $d(m,n)=\log_2 f(m,n+1) - \log_2 f(m,n)$.
Then Conjecture~\conjConv\ implies $d(m,n)\le d(m,n+1)\le d(m,n+2)\le\dots$
whence $\log_2 f(m,n+k)-\log_2 f(m,n) \ge k d(m,n)$. Dividing by $km$ and
passing to the limit
when $k\to\infty$, we obtain $c_m\ge d(m,n)/m=(n+1)c(m,n+1)-nc(m,n)$.
\qed
\enddemo


\head\sectCtwo. The exact value of $c_2$ (proof of Theorem 1)
\endhead
For $a,c\ge 0$, $a\equiv c\mod 2$,
let $g^*_{a,c}$ be the number of primitive lattice triangulations
of the trapezoid $T(a,c)$ spanned by $(0,0)$, $(a,0)$, $(1,2)$, $(1+c,2)$
(if $a=0$ or $c=0$, then $T(a,c)$ degenerates to a triangle).
When $a\not\equiv c\mod 2$, we set $g^*_{a,c}=0$.
We also set $g^*_{0,0}=1$.
Let $G^*(x,z)$ be the generating function for $g^*_{a,c}$:
$$
\split
     G^*(x,z) &= \sum_{a,c\ge0} g^*_{a,c}\, x^a z^c\\
     &= 1 + (x^2 + xz + z^2) + (6x^4 + 10x^3 z + 12 x^2 z^2 + 10 x z^3 + 6z^4) + \dots
\endsplit
$$
Let $g^*_n$ be the coefficient of $x^{2n}$ in the series
$G^*(x,x)=\sum_{n\ge0} g^*_n\, x^{2n}$, i.e.
$$
    g^*_n = g^*_{0,2n} + g^*_{1,2n-1} + g^*_{2,2n-2} + \dots + g^*_{2n,0}.
$$
Then Theorem 1 follows immediately from Lemmas \lemFone\ and \lemFtwo\ below.

\proclaim{ Lemma \lemFone }
 $\lim_{n\to\infty} f(2,n)^{1/n} = \lim_{n\to\infty} (g_n^*)^{1/n}$.
\endproclaim

\demo{ Proof } The rectangle $2\times(n-1)$ can be placed into $T(n,n)$, hence
$f(2,n-1) < g^*_n$. On the other hand, the union of $T(a,c)$ with its image under
the central symmetry with center $(\tfrac12(a+c+1),1)$ is $T(a+c,a+c)$, and it can be placed into
the rectangle $2\times(a+c+1)$, hence $(g^*_{a,c})^2 < f(2,a+c+1)$. Therefore
$$
    \frac{g^*_n}{2n}=\!\sum_{a+c=2n}\!\!\!\!\frac{g^*_{a,c}}{2n}
    \le\max_{a+c=2n}g^*_{a,c}\le f(2,2n+1)^{1/2}\le(g^*_{2n+2})^{1/2}
$$
whence $\frac1n\big(\log g^*_n - \log(2n)\big) \le \tfrac1{2n} f(2,2n+1)
    \le \tfrac1{2n} g^*_{2n+2}$
and the result follows because $\frac1n\log(2n)\to 0$.
\qed\enddemo

\proclaim{ Lemma \lemFtwo } $\lim_{n\to\infty}(g^*_n)^{1/n} = \alpha$ where $\alpha$ is as in Theorem 1.
\endproclaim

\demo{ Proof }
For $a,c\ge 0$, $a\equiv c\mod 2$,
let $g_{a,c}$ be the number of those primitive lattice triangulations
of the trapezoid $T(a,c)$ which do not contain interior edges of the form
$[(k,0),(l,2)]$, in other words, primitive lattice triangulations which 
agree with the subdivision of $T(a,c)$ into two triangles and two trapezoids
depicted in Figure \figT(left). If $a+c$ is odd, we set $g_{a,c}=0$.
By convention, we set $g_{0,0}=0$.
Let $G(x,z)=\sum_{a,c\ge 0} g_{a,c}\,x^a z^c$ be the generating function.

\midinsert
\centerline{\epsfxsize=100mm\epsfbox{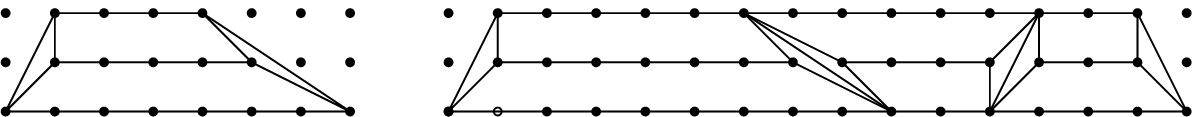}}
\botcaption{Figure \figT}
\endcaption
\endinsert

The edges of the form $[(k,0),(l,2)]$ of
any primitive lattice triangulation cut $T(a,c)$ into smaller trapezoids. They
can be transformed into $T(a_i,c_i)$'s with $\sum a_i = a$ and $\sum c_i = c$
by uniquely determined lattice automorphisms of the form $(x,y)\mapsto(x+p_iy+q_i,y)$
with $p_i,q_i\in\Z$
(see Figure \figT). Hence
$$
     g^*_{a,c} =
     \sum_{\smallmatrix a_1+\dots+a_k=a\\c_1+\dots+c_k=c\endsmallmatrix}\;
          \prod_{j=1}^k g_{a_j,c_j},
\quad\text{thus}\quad
     G^*(x,z) = \frac1{1 - G(x,z)}.                           \eqno(\eqF)
$$

Easy to see (cf.~(\eqCone)) that the number of primitive lattice triangulations of the
narrow (i.e. of width 1) trapezoids in Figure~\figT\ are
binomial coefficients, hence
$G(x,z)=(x^2+xz+z^2) + (5x^4 + 8x^3z + 9x^2z^2 + 8xz^3 + 5z^4)+\dots$.

One can also check that $g_{a,c}=f_{a,(a+c)/2-1,c}$ where $f_{a,b,c}=f(S_{a,b,c})$
are the numbers discussed in Example~\exTetrisA.
Hence (cf.~Example~\exTetrisB)
$$
\split
  G(x,z)&=\sum_{a,c} f_{a,(a+c)/2-1,c}\,x^az^c
        =\text{coef}_{u^0}\Big[\sum_{a,b,c}\,f_{a,b,c}\,x^a u^{2b-a-c+2}z^c\Big]
\\&
     =\text{coef}_{u^0}\Big[u^2\sum_{a,b,c}f_{a,b,c}
            \,(x/u)^a u^{2b}(z/u)^c\Big]
     =\text{coef}_{u^{-1}}\Big[u F\big(x/u,u^2,z/u)\Big].
\endsplit  
$$
Since the function $1/(1-x-y-z+xz)=1/\big((1-x)(1-z)-y\big)$
is analytic in the domain $\max\big(|x|,|y|,|z|\big)<1/2$, its
power series $\sum f_{a,b,c}\,x^ay^bz^c$ (see Example~\exTetrisA) converges to it
in this domain. Therefore, for $0<\varepsilon\ll r<1/2$,
 the Laurent series of $F(x/u,u^2,z/u)$ converges in the domain
$\max\big(|x|,|z|\big)<\varepsilon$, $r-\varepsilon<|u|<r+\varepsilon$.
Hence, for $x$ small enough, we have
$$
   G(x,x) = \text{coef}_{u^{-1}}\big[F(x/u,u^2,x/u)\big]
  =\frac{1}{2\pi i}\oint_{|u|=r}\frac{u\,du}{(1-x/u)^2-u^2}
$$
and
$$
\split
    \frac{u}{(1-x/u)^2-u^2} &=  -\frac{u}{2(u^2+u-x)} - \frac{u}{2(u^2-u+x)}
\\&
 = \sum_{j=1}^2\frac 1{2(u_j^+-u_j^-)}\Bigg(\frac{u_j^+}{u-u_j^+}+\frac{u_j^-}{u-u_j^-}\Bigg),
\endsplit
$$
where, for $|x|$ small enough,
$$
   u_1^\pm = -\tfrac12(1\pm\sqrt{1+4x}),\quad
   u_2^\pm =  \tfrac12(1\pm\sqrt{1-4x});\qquad |u_j^+|>r,\; |u_j^-|<r.
$$
Thus
$$
   G(x,x) = \sum_{j=1}^2\!\underset{\;\;u=u_j^-}\to{\text{Res}}\!\Big(\dots\Big)
          = \sum_{j=1}^2\frac{u_j^-}{2(u_j^+ - u_j^-)}
            = \frac1{4\sqrt{1-4x}} + \frac1{4\sqrt{1+4x}} - \frac12.
$$
The graph of the function $y = G(x,x)$ sits in the algebraic curve
$$
  (2y+1)^2(16x^2-1)\big(4x^2 + (y^2+y)(16x^2-1)\big)+x^2 = 0.
$$
By (\eqF), the poles of $G^*(x,x)$ are the $x$-coordinates of the intersections of this curve
with the line $y=1$, i.e., the roots of $5184x^4 - 611 x^2 + 18$ (the smallest ones being $\pm\sqrt{1/\alpha}$),
and the branching points are $\pm1/4$. Hence the radius of convergence of the series
$G^*(x,x)=\sum g^*_n x^{2n}$ is $\sqrt{1/\alpha}$
whence $\lim_{n\to\infty} (g^*_n)^{1/n} = \alpha$.
\qed\enddemo


\head\sectCthree. Computation of $c_3$ (proof of Theorem 2)
\endhead

\subhead\sectCthree.1. Preparation
\endsubhead
For $a,d\ge0$ such that $a\not\equiv d+1\mod3$, let
$h^*_{a,d}$ be the number of primitive lattice triangulations of the
trapezoid $T_3(a,d)$ spanned by $(0,0),(1,3),(1+d,3),(a,3)$. We set $h^*_{0,0}=1$
and $h^*_{a,d}=0$ when $a\equiv d+1\mod 3$, and
we consider the generating function 
$$
   H^*(x)=\sum_n h^*_n x^n
         = \sum_{a,d\ge 0}h^*_{a,d}\, x^{a+d} = 1 + x + 3x^2 + 19x^3 + 125 x^4 +\dots
$$
Similarly to the beginning of proof of Lemma~\lemFtwo,
we define $h_{a,d}$ as the number of the triangulations of $T_3(a,d)$
which do not have edges of the form $[(k,0),(l,3)]$ and we consider the generating function
$$
   H(x)=\sum_n h_n x^n = \sum_{a,b\ge 0}h_{a,d}\, x^{a+d}\,
                = x + 2x^2 + 14 x^3 + 86 x^4 + 712 x^5 +\dots
$$
These functions satisfy the relation similar to (\eqF) specialized for $x=z$: 
$$
    H^*(x) = 1/(1-H(x))
$$
Indeed, the edges of the form $[(k,0),(l,3)]$ cut $T_3(a,d)$ into smaller trapezoids.
Each of them can be mapped to a standard one by a unique lattice automorphism of the form 
$(x,y)\mapsto(x+py+q,y$ or $3-y)$ with $p,q\in\Z$
(in contrary to \S2, here the upper and lower horizontal sides of the trapezoids are mixed,
so we do not have (\eqF) for two-variable generating functions).
In Figure \figTT\ we illustrate the relation
$$
     h^*_3 = h^*_{03} + h^*_{12} + h^*_{30}
      = h_{01}^3 + 2h_{01}(h_{11}+h_{20}) +  (h_{03} + h_{12} + h_{30})
      = h_1^3 + 2 h_1 h_2 + h_3.
$$

\midinsert
\centerline{\epsfxsize=112mm\epsfbox{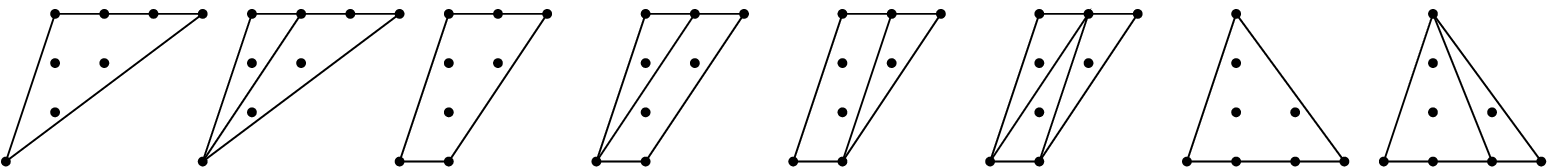}}
\centerline{
$h_{03}\,\;\;+\;\;h_{01}h_{20}\,+\,h_{12}\;\;+\;\;h_{01}h_{11}+h_{11}h_{01}
     \;+\;h_{01}^3\;\;\;+\;\;\;h_{30}\;\;+\;\;h_{20}h_{01}$ }
\botcaption{ Figure \figTT }
\endcaption
\endinsert
Similarly to Lemma \lemFone, we have
$\lim_n f(3,n)^{1/n} = \lim_n(h^*_{2n})^{1/n} = 1/\beta^{2}$
where $\beta$ is the real positive root of the equation $G(x)=1$,
hence $c_3=-\frac23\log\beta$.

\subhead\sectCthree.2. Recurrent relations
\endsubhead
Using the notation introduced in \S\sectRecursion,
let us set $m=3$, $\varphi_0(x)=\frac13x-1$, and
$$
   F(x,y,z,w)=\sum_{a,b,c,d} f_{a,b,c,d}\,x^ay^bz^cw^d,
$$
$$
   G_1(x,z,w)=\sum_{a,c,d} g^{(1)}_{a,c,d}\,x^a z^c w^d\,,
\qquad
   G_2(x,y,w)=\sum_{a,b,d} g^{(2)}_{a,b,d}\,x^a y^b w^d
$$
$$
   H_k(x,w)=\sum_{a,d}g^{(k)}_{a,d}\,x^a w^d\,,
\qquad\quad (k=1,2)
$$
where all the coefficients are of the form $f(S)$ (see \S\sectRecursion) for the
$\varphi_0$-admissible shapes in Figure~\figFGH\
where $(0,a)$, $(1,b)$, $(2,c)$, and $(3,d)$ (if present) are the coordinates
of integral points on the upper part of the boundary of $S$.
The lower corners of $S$ are at the points $(0,-1)$ and $(3,0)$.
If the congruences
given in Figure~\figFGH\ are not satisfied, then the corresponding numbers are zero.
If $\min(a+1,b,c,d)<0$, they are also zero (this case does not correspond to any
$\varphi_0$-admissible shape). By convention, we also set $h_{-1,0}^{(2)}=0$
(the case when $S$ degenerates to a segment).

\midinsert
\centerline{\epsfxsize=100mm\epsfbox{fgh.eps}}
\centerline{
 $f_{a,b,c,d}$  \hskip 10mm
  $g^{(1)}_{a,c,d}$ \hskip 12mm
  $g^{(2)}_{a,b,d}$ \hskip 13mm
  $h^{(1)}_{a,d}$ \hskip 14mm
  $h^{(2)}_{a,d}$ \hskip 1mm
}
\centerline{%
  \hskip 25mm $c-a\equiv 1(2)$
  \hskip 2mm $d-b\equiv 1(2)$
  \hskip 2mm $d-a\equiv 2(3)$
  \hskip 2mm $d-a\equiv 1(3)$
}
\botcaption{Figure \figFGH}
\endcaption
\endinsert

In terms of the generating functions, the recurrent relations in
Lemma~\lemTetrisB\ read (cf.~Examples~\exTetrisA, \exTetrisB):
$$
\split
   F(x,y,z,w)&(1-x-y-z-w+xz+yw+xw)
\\
  &=y^{1/2}(1-w)G_1(xy^{1/2},y^{1/2}z,w)
  +z^{1/2}(1-x)G_2(x,yz^{1/2},z^{1/2}w),
\endsplit
$$
$$
\split
  &G_1(x,z,w)(1-w)=\text{coef}_{u^{-1}}\big[F(x/u,u^2,z/u,w)(1-w)\big]+x^{-1},
\\
  &G_2(x,y,w)(1-x)=\text{coef}_{u^{-1}}\big[F(x,y/u,u^2,w/u)(1-x)\big]
\endsplit
$$
(the asymmetry between $G_1$ and $G_2$ is caused by the asymmetry of $\varphi_0$),
$$
\split
  &H_1(x,w)=\text{coef}_{u^{-1}}\big[ G_1(x/u,u^3,w/u^2)\big],\\
  &H_2(x,w)=\text{coef}_{u^{-1}}\big[ G_2(x/u^2,u^3,w/u)\big].
\endsplit
$$
Notice that in this subsection, by generating functions we mean formal series.
Let us consider the symmetrized generating functions
$$
\split
 \tilde F(x,y,z,w) &= F(x,y,z,w)+F(w,z,y,x),\\
 \tilde G(x,z,w)   &= G_1(x,z,w) + G_2(w,z,x),\\
 \tilde H(x,w) &= H_1(x,w)+H_2(w,x).
\endsplit
$$
The above relations for $F,G_1,G_2,H_1,H_2$
imply immediately:
$$
\matrix
   \tilde F(x,y,z,w)
    ( 1-x-y-z-w+xz+yw+xw )
  \qquad\qquad\qquad\qquad\qquad\qquad
\\
  \text{\hbox{\vbox to 14pt{}}}
  =   y^{1/2}(1-w)\tilde G(xy^{1/2},y^{1/2}z,w)
     +z^{1/2}(1-x)\tilde G(x,yz^{1/2},z^{1/2}w),
\endmatrix
   \eqno(\eqFF)
$$
$$
  \tilde G(x,z,w)(1-w)=\text{coef}_{u^{-1}}\big[\tilde F(x/u,u^2,z/u,w)(1-w)\big]+x^{-1}
    \eqno(\eqGG)
$$
$$
  \tilde H(x,w)=\text{coef}_{u^{-1}}\big[\tilde G(x/u,u^3,w/u^2)\big].\eqno(\eqHH)
$$


\subhead\sectCthree.3. The equation
\endsubhead
We are going to obtain an equation for $\tilde G(xt^{-1/2},t^{3/2},t^{-1}x)$ by
expressing $\tilde F$ via $\tilde G$ from (\eqFF) and plugging it
to (\eqGG). To this end we need to divide power series by polynomials.
 However, when some variables appear with powers varying from $-\infty$ to $+\infty$,
the meaning of such division should be precised.
To illustrate a possible ambiguity, let us consider the
expression $\text{coef}_{u^{-1}}\big[1/(x-uy)\big]$. It can be understood either as
$$
  \text{coef}_{u^{-1}}\Big[ \frac{x^{-1}}{1-uyx^{-1}}\Big]
  =\frac1x\text{coef}_{u^{-1}}\Big[1+\frac{uy}x+\frac{u^2y^2}{x^2}+\dots\Big]
  =0
$$
or as
$$
  \text{coef}_{u^{-1}}\Big[ -\frac{(uy)^{-1}}{1-x(uy)^{-1}}\Big]
  =-\text{coef}_{u^{-1}}\Big[\frac1{uy}\Big(1+\frac{x}{uy}+\frac{x^2}{u^2y^2}
    +\dots\Big)\Big]
  =-\frac1y.
$$
To avoid this kind of ambiguity, we introduce a new formal variable $q$
and consider the formal series
$$
\split
   F_q(x,y,z,w)&=F(xq,yq^2,zq^2,wq),\\
   G_{1,q}(x,z,w)&=G_1(xq^2,zq^3,wq),\\
   G_{2,q}(x,y,w)&=G_2(xq,yq^3,wq^2),\\
   H_{k,q}(x,w)&=H_k(xq^3,wq^3),\qquad k=1,2,
\endsplit
$$
and all the generating functions will be treated as elements of
the ring
$$
   \Bbb Z[x^{\pm1},y^{\pm1/2},z^{\pm1/2},w^{\pm1},u^{\pm1/2},t^{\pm1/2}]((q))
$$
of formal power series in $q$ (starting, maybe, with a negative power) whose
coefficients are Laurent polynomials in $x,y^{1/2},\dots$.

The geometric meaning of an exponent of $q$ is twice the doubled signed area of the
$\varphi_0$-admissible shape corresponding a monomial, i.e.,
$2\int_0^3\varphi(x)\,dx$ where the graph of $\varphi$ is the upper boundary
of the shape.
One can easily check by hand that
$$
\split
  F_q &=  (xq)^{-1} + (1+x^{-1}w) + (x+w+x^{-1}y + x^{-1}z + x^{-1}w^2)q+\dots
\\
  G_{1,q} &= x^{-1}q^{-2} + w(xq)^{-1} + w^2x^{-1} + w^3x^{-1}q
     + (x+w^4x^{-1})q^2+\dots
\\
  G_{2,q} &= x^{-1}wq + (w+yx^{-1})q^2 + (wx+2y)q^3 + (wx^2+4xy)q^4+\dots
\\
  H_{1,q} &=  wx^{-1} + xq^3 + 4w^2 q^6
  +(30wx^2 + 24 w^4 x^{-1})q^9
 +\dots
\\
  H_{2,q} &= wq^3 + 5(x^2+w^3x^{-1})q^6+ 32w^2x\,q^9 + \dots
\endsplit
$$

Further, we define $\tilde F_q$, $\tilde G_q$, $\tilde H_q$
by the same formulas as in \S\sectCthree.2 but with the subscript $q$
everywhere. For example,
$$
 \tilde G_q(x,z,w)=\frac{1}{xq^2} + \frac w{xq} + \frac{w^2}x
    + \Big(\frac{w^3}x + \frac xw\Big)q
     + \Big(2x+\frac{w^4}x + \frac zw\Big)q^2+\dots
$$
Then the relations (\eqFF)--(\eqHH) take the form
$$
   \tilde F_q
  =\frac{
      qy^{1/2}(1-wq)\tilde G_q(xy^{1/2},y^{1/2}z,w)
     +qz^{1/2}(1-xq)\tilde G_q(wz^{1/2},z^{1/2}y,x) }
    { 1-xq-yq^2-zq^2-wq+xzq^3+ywq^3+xwq^2 },
   \eqno(\eqFFq)
$$
$$
  \tilde G_q(x,z,w)=\text{coef}_{u^{-1}}\big[q\tilde F_q(x/u,u^2,z/u,w)\big]
     +\frac1{x(1-wq)q^2},
    \eqno(\eqGGq)
$$
$$
  \tilde H_q(x,w)=\text{coef}_{u^{-1}}\big[q\tilde G_q(x/u,u^3,w/u^2)\big].\eqno(\eqHHq)
$$

Let us set
$$
\split
      g_q(x,t) &= t^{1/2}x^2q^2 \tilde G_q(x^2t^{-1/2},x^3t^{3/2},xt^{-1})
    \\&  = t + xq + t^{-1}x^2q^2 + (t^{-2}+t)x^3q^3
       +(t^{-3}+2+t^3)x^4 q^4
      +\dots
\endsplit
$$
The parity condition on the indices of nonzero coefficients of $G_1$ and $G_2$
(see Figure~\figFGH)
ensures that the series $g_q(x,t)$ does not have fractional powers.
Moreover, $x$ and $q$ appear in each monomial of $g_q$ with the same power, thus we have
$g_q(x,t)=g(xq,t)$ with $g(x,t)\in\Z[t^{\pm1}]((x))$. 

By plugging (\eqFF) into (\eqGG), denoting the denominator in (\eqFFq) by
$Q_q(x,y,z,w)$, and observing that
$$
  \text{coef}_{u^{-1}}\big[\Cal F(x,t,u)\big]
    =\text{coef}_{u^{-1}}\big[xt^{-1/2}\Cal F(x,t,uxt^{-1/2})\big]
     \eqno(\eqPW)
$$
for any formal Laurent series in $u$,
we obtain
$$
\split
  g_q(x,t) 
      &\overset\text{(\eqGGq)}\to=
      \text{coef}_{u^{-1}}\Big[t^{1/2} x^2 q^3
      \tilde F_q\Big(\frac{x^2}{ut^{1/2}},u^2,\frac{x^3t^{3/2}}{u},\frac{x}{t}\Big)\Big]
         +\frac{t^2}{t-xq}
\\
      &\overset\text{(\eqPW)}\to=
      \text{coef}_{u^{-1}}\Big[x^3 q^3
      \tilde F_q\Big(\frac{x}{u},\frac{x^2u^2}{t},\frac{x^2t^2}{u},\frac{x}{t}\Big)\Big]
         +\frac{t^2}{t-xq}
\\
  &\overset\text{(\eqFFq)}\to=
     x^2q^2\text{coef}_{u^{-1}}\Bigg[
       \frac{  \frac{u}{t}\big(1-\frac{xq}{t}\big)g_q(x,t)
             + \frac{t}{u}\big(1-\frac{xq}{u}\big)g_q(x,u) }
      {Q_q\big(x/u,\,x^2u^2/t,\,x^2t^2/u,\,x/t\big)} \Bigg]
         +\frac{t^2}{t-xq}
\\
  &= x^2q^2\text{coef}_{u^{-1}}\Big[
       \frac{  u^3(t-xq)g_q(x,t)
             + t^3(u-xq)g_q(x,u) }
      {P(xq,t,u)} \Big]
         +\frac{t^2}{t-xq}
\endsplit
$$
where
$$
   P(x,t,u) = u^2t^2-(u+t)utx +(1-t^3-u^3)utx^2+ (t^4+u^4)x^3.
  \eqno(\eqDefP)
$$

We see that the variables $x$ and $q$ are ``synchronized''
in the right hand side of the obtained equation:
they occur with the same power in each monomial
of each power series in this expression.
Hence we obtain the following identity in the ring $\Z[t^{\pm1},u^{\pm1}]((x))$:
$$
    g(x,t)\Psi(x,t)
    = 
     \frac{t^2}{t-x}
      + \text{coef}_{u^{-1}}\left[
        \frac{t^3x^2(u-x)\,g(x,u)}{P(x,t,u)}\right]
    \eqno(\eqg)
$$ 
where
$$
  \Psi(x,t) = 1-x^2(t-x)\Phi(x,t),
\qquad
    \Phi(x,t)=\text{coef}_{u^{-1}}\big[u^3/P(x,t,u)\big].
$$
Here are several initial terms of these series:
\footnote{
  All coefficients of $\Phi$ and $1-\Psi$ that I have computed
  are positive. If they are really all positive, it would be interesting
  to find their combinatorial meaning.
}
$$
  \Phi(x,t)=t^{-2}x^2+(t^{-3}+1)x^3+(t^{-4}+2t^{-1}+t^2)x^4
            +(t^{-5}+3t^{-2}+3t)x^6+\dots
$$
$$
  \Psi(x,t)=1-t^{-1}x^4-tx^5-(1+t^3)x^6-(t^{-1}+2t^2)x^7-(6t^{-2}+3t)x^8-\dots
  \eqno(\eqPsi)
$$

Having found $g$ form (\eqg), we can compute $\tilde H(x,x)$. Indeed, by (\eqHH) we have
$$
   x\tilde H_q(x^3,x^3)
    =\text{coef}_{t^0}\big[tx\tilde G_q(x^3/t,t^3,x^3/t^2)\big]
    =\text{coef}_{t^0}\big[t^{1/2}x\tilde G_q(x^3/t^{1/2},t^{3/2},x^3/t)\big].
$$
Replacing $t$ by $x^2t$ (cf.~(\eqPW)) and setting $q=1$, we obtain
$$
        x \tilde H(x^3,x^3) = \text{coef}_{t^{0}}\big[g(x,t)\big].
     \eqno(\eqHg)
$$


\subhead\sectCthree.4. Computation
\endsubhead
In this subsection we study the analytic functions defined by the series
discussed in the previous subsection.

By \S\sectCthree.1, we need to find the smallest positive
pole of $H^*(x)$, that is the smallest positive zero $\beta$ of $1-H(x)$.
One can check that
$$
      H(x) = x\tilde H(x,x).
    \eqno(\eqHs)
$$
Being the sum of a power series with positive coefficients,
the function $x\tilde H(x,x)$ is increasing when $x>0$,
thus it is enough know how to compute with any given precision
the value of $\tilde H(x,x)$ for any fixed $x$
in an interval containing $\beta$.
By (\eqHg), this can be done by numerical integration of the function $g(x^{1/3},t)$
along a suitable contour $\Gamma_x$ (cf.~the proof of Lemma~\lemFtwo).
Thus we need to be able to compute $g(x,t)$ for any $x\in[0,x_0^+]$ and
$t\in\Gamma_x$ for some $x_0^+ > x_0=\beta^{1/3}$.
This can be done because for a fixed $x$, after replacing $\text{coef}_{u^{-1}}[\dots]$
by $\frac1{2\pi i}\int_{\Gamma_x}(\dots)du$,
the equation (\eqg)
becomes a Fredholm equation for the function $g$ restricted to $\Gamma_x$.
Now we pass to more detailed explanations.

Let
$$
\split
   &\Gamma=\{(x,t,u)\in\R\times\C^2\mid 0<x<1/2,\; |t|=|u|=1\},
\\
   &\Gamma'=\{(x,t)\in\R\times\C\mid 0<x<1/2,\; |t|=1\}.
\endsplit
$$

\proclaim{ Lemma \lemP }
The polynomial $P(x,t,u)$ defined in (\eqDefP) does not vanish on $\Gamma$.
For any fixed $(x,t)\in\Gamma'$,
the polynomial $P(x,t,u)$ has two simple roots $u_k(x,t)$, $k=1,2$,
in the unit disk $|u|<1$ and two simple roots outside it.
\endproclaim

\demo{Proof}
The first statements can be checked
using any software for symbolic computations.
This can be done, for example, as follows.
Let $S^1$ be the unit circle in $\C$.
Then $\Gamma=(0,1/2)\times S^1\times S^1$.
We can identify $S^1$ with $\RP^1$
by some rational parametrization.
Then $\Re P$ and $\Im P$ become real rational functions on the variety $\Gamma$ and,
by computing resultants, discriminant, etc., one can check that
the real algebraic curve given by the equations $\Re P=\Im P=0$
does not enter in the layer $0<r<1/2$. More precisely, let
$p(x,T,U)$ and $q(x,T,U)$ be real polynomials such that
$$
   P(x,\zeta(T),\zeta(U))=\frac{p(x,T,U)+iq(x,T,U)}{(i+T)^4(i+U)^4},
\qquad
 \zeta(X)=\frac{i-X}{i+X}.
$$
Note that $\zeta(\R)=S^1\setminus\{-1\}$, hence $(x,T,U)$ are coordinates
on the affine chart $\Gamma\setminus\{(t+1)(u+1)=0\}$ of $\Gamma$.
The projection of the real algebraic curve $\Gamma\cap\{P=0\}$ onto the plane $(x,T)$
is given by the equation $R(x,T)=0$
where $R(x,T)$ is the resultant of $p$ and $q$ with respect to $U$.
To prove that the curve $R(x,T)=0$ does not have real points with $0<x<1/2$,
we compute the real roots of $D(x)=0$ on this interval where $D(x)$ is the
discriminant of $R$ with respect to $T$,
and we check that the equations $R(x_k,T)=0$ for each $k=1,\dots,2n+1$
do not have real roots where $0<x_1<\dots<x_{2n+1}<1/2$ and 
$x_k$ with even $k$ are all the real roots of $D(x)$ on the interval $0<x<1/2$.
This computation shows that
$P(x,u,t)\ne0$ when $(x,u,t)\in\Gamma$ and $(t+1)(u+1)\ne0$.
Then we check that $P(x,\zeta(T),-1)\ne0$,
$P(x,-1,\zeta(U))\ne0$, and $P(x,-1,-1)\ne0$ for $0<x<1/2$, $T\in\R$.

Similarly one can check that for any fixed $(x,t)\in\Gamma'$,
the discriminant of $P$ with respect to the variable $u$
does not vanish, hence for any fixed $(x,t)\in\Gamma'$,
all the four roots of $P$ (viewed as a polynomial
in $u$) are pairwise distinct.

Therefore, the number of roots of $P$
in the unit disk $|u|<1$ is constant.
Thus, to prove the second statement, it is enough to check
it for some value of $x$ and $t$, for example, for $t=1$
and a very small $x$.
\qed\enddemo

\proclaim{ Lemma \lemPhi }
(a).
The formal power series $1/P(x,t,u)\in\Z[t^{\pm1},u^{\pm1}]((x))$
converges to the function $1/P(x,t,u)$ in a neighborhood
of $\Gamma\cap\{x<\frac{1}{4}\}$.

\smallskip
(b). The formal power series $\Phi(x,t)\in\Z[t^{\pm1}]((x))$ converges to an
analytic function (which we also denote by $\Phi(x,t)$) in a neighborhood of
$\Gamma'\cap\{x<\frac{1}{4}\}$.
The function $\Phi(x,t)$ admits an analytic continuation
to a neighborhood of $\Gamma'$ defined by the Cauchy integral
$$
  \Phi(x,t) = \frac{1}{2\pi i}\oint_{|u|=1}\frac{u^3\,du}{P(x,t,u)}
  =\sum_{k=1}^2\frac{u_k(x,t)^3}{P'_u(x,t,u_k(x,t))}
  \eqno(\eqPhi)
$$
where $u_1(x,t)$ and $u_2(x,t)$ are the roots of $P$ in the unit disk $|u|<1$;
see Lemma~\lemP.
\endproclaim

\demo{ Proof } The power series
$1/P(x,t,u)$ involved in the definition $\Phi(x,t)$ is a power series expansion
with respect to $x$, hence
$1/P=a_0^{-1}(1+X+X^2+\dots)$ where $X=(a_1+a_2+a_3)/a_0$ and
$a_k=x^k\text{coef}_{x^k}[P]$. If $(x,t,u)\in\Gamma$, then
$|a_0|=1$, $|a_1|\le 2x$, $|a_2|\le 3x^2$, $|a_3|\le 2x^3$, and
thus $|X|\le 2x+3x^2+2x^3$. Therefore $|X|<1$ for $x<1/4$,
whence the convergence of $1/P$ in the required domain.
This fact combined with Lemma~\lemP\ implies all the other assertions
of the lemma.
\qed\enddemo

Mathematica function {\,\tt Psi\,} in Figure~\figProg\ computes
$\Psi(x,t)$ for $(x,t)\in\Gamma'$ with any given precision.

Notice that one of the functions $u_1(x,t)$ or $u_2(x,t)$
has a ramification point at $(x,t)=(1/2,1)$, and hence the functions
$\Phi$ and $\Psi$ are ramified in this point as well.
The Laurent-Puiseux expansion of $\Psi(x,1)$ in powers of $s=\sqrt{1/2-x}$ is
$$
   \Psi(x,1)=-\tfrac{1}{4\sqrt6}s^{-1}
  +\tfrac{12-\sqrt2}{8}
  -\tfrac{3}{8\sqrt6}s
  -\tfrac{3}{8\sqrt2}s^2 + \tfrac{103}{96\sqrt6}s^3 - \tfrac{87}{32\sqrt2}s^4
  +\tfrac{2635}{192\sqrt6}s^5
   +\dots
$$

Let $x_0^-=\frac{16}{33}$
and $x_0^+=\frac{17}{35}$. We shall see later that
$x_0\in[x_0^-,x_0^+]$; in fact, $x_0^{\pm}$ are given by initial
segments of the continued fraction of $x_0$.

Using the expansion of $\Psi$ at $(\tfrac12,1)$ and
computing the values of
$\Psi(x,t)$ (with the program in Figure~\figProg)
on a sufficiently dense grid on $\Gamma'$, one can check that
$\Psi$ does not vanish on $\Gamma'\cap\{x<x_0^+\}$ 
and
$$
    \min_{0\le x\le x_0^+,\,|t|=1}|\Psi(x,t)|
    =\min_{0\le x\le x_0^+,\,|t|=1}\Re\Psi(x,t)
    =\Psi(x_0^+,1)=0.44768...
                                  \eqno(\eqMinPsi)
$$
See the level lines of $\Re\Psi$ in Figure~\figPsiLevel;
we omit the details of the error estimate.

Using Lemma~\lemTaylor\ applied to the function
$|P(x/4,e^{i\tau},e^{i\theta})|^2$
with an appropriately chosen $h$, we find
$$
   \min_{x<x_0^+,\,|t|=|u|=1}|P| = P(x_0^+,1,1) = 0.02183...  \eqno(\eqMinP)
$$
(here we rescaled $x$ to equilibrate the partial derivatives).
The computation can be fastened by choosing different grid
in different zones of $\Gamma$.
In our computation, the grid step
varied from $h=1/300$ near the
point of minimum to $h=1/20$ far from it. To estimate the error, we
used evident coarse bounds for the fourth derivatives and, using them,
computed finer upper bounds for the second derivatives in each zone
again using Lemma~\lemTaylor.

\midinsert
\centerline{\epsfxsize=120mm\epsfbox{psi-re.eps}}
\botcaption{Figure \figPsiLevel}
   Level lines of $\Re\Psi(x,t)$  for $|t|=1$. The
   shown vertical line is $x=x_0$ or $x=x_0^+$ (no difference with
this resolution).
\endcaption
\endinsert

\proclaim{ Lemma \lemG }
The formal series $g(x,t)$ (introduced in \S\sectCthree.2)
converges in some neighborhood of $\Gamma'\cap\{|x|<2^{-3/2}\}$.
\endproclaim

\demo{ Proof }
By Anclin's theorem [\refAnclin], the number of primitive lattice triangulations
of a lattice polygon $\Pi$ is bounded above by
$2^N$ for $N=\#\big(\Pi\cap(\Z^2\setminus\frac12\Z^2)\big)$, and it is easy to derive
from Pick's formula that $N<3\text{Area}(\Pi)-3/2$.
The area of the shape corresponding to $g^{(k)}_{a,c,d}$ is $(2a+3c+d+3)/2$.
Hence $\tilde g_{a,c,d} < c_0 2^{3(2a+3c+d)/2}$ for some constant $c_0$ and,
for $|t|=1$, we obtain
$$
\split
   |g(x,t)|&\le
   x^2\sum_{a,c,d} \big|\tilde g_{a,c,d}\, x^{2a} x^{3c} x^d\big|
\\&
   \le c_0x^2\sum_{a,c,d} \big|2^{3(2a+3c+d)/2}\,
         x^{2a+3c+d}\big|
    =c_0x^2\sum_n 2^{3/2n}A_n x^n,
\endsplit
$$
where $A_n=\#\{(a,c,d)\in\Z_+^3\mid 2a+3c+d=n\}$. Since $A_n$ is bounded by
a polynomial function of $n$, the series converges for $x<2^{-3/2}$.
\qed
\enddemo

Lemmas \lemPhi\ and \lemG\ combined with (\eqg) and (\eqMinPsi)
imply that the function $g(x,t)$
is analytic in a neighborhood of $\Gamma'\cap\{x<2^{-3/2}\}$, and it satisfies
the condition
$$
    g(x,t) = \frac{t^2}{(t-x)\Psi(x,t)}
        +\frac{1}{2\pi i}\oint_{|u|=1}\frac{x^2t^3(u-x)g(x,u)\,du}{P(x,t,u)\Psi(x,t)}.
   \eqno(\eqFredG)
$$
For any fixed $x$, this is a Fredholm equation of the second kind 
for $g(x,t)$ considered as a function of $t$.

\proclaim{ Lemma \lemFredG }
The function $g(x,t)$ analytically extends to a neighborhood of
$\Gamma'\cap\{x<x_0^+\}$ and it satisfies the equation (\eqFredG) in this domain.
\endproclaim

\demo{ Proof } Let us rewrite (\eqFredG) in a more conventional form
$$
   \varphi_g(x,\tau)=f(x,\tau)
     +\int_0^1 K(x,\tau,\theta)\varphi_g(x,\theta)\,d\theta     \eqno(\eqFredTau)
$$
where we set $t=e^{2\pi i\tau}$, $u=e^{2\pi i\theta}$, and
$$
   \varphi_g(x,\tau)=g(x,t),\quad f(x,\tau)=\frac{t^2}{(t-x)\Psi(x,t)},\quad
  K(x,\tau,\theta)=\frac{x^2t^3u(u-x)}{P(x,t,u)\Psi(x,t)}.
$$

As we already pointed out, $g$ satisfies (\eqFredG) and thus $\varphi_g$
satisfies (\eqFredTau) for small $x$.
Thus, thanks to the Identity Theorem for analytic functions,
it is enough to show that for any $x\in[0,x_0^+]$ there exists
a unique solution of (\eqFredTau)
and that it is analytic with respect to $(x,\tau)$.
Hence, by Lemma~\lemParam, it suffices to show that $1$ is not an eigenvalue
of $\Cal K_x$ for any $x\in[0,x_0^+]$ where
$\Cal K_x:\Cal C[0,1]\to\Cal C[0,1]$ is the the Fredholm integral operator
which takes $\varphi(\tau)$ to
$\psi(\tau)=\int_0^1 K(x,\tau,\theta)\varphi(\theta)\,d\theta$.
The latter fact, in its turn, follows from the bound
$$
   \max_{0\le x\le x_0^+}\Cal N_2(x)=\Cal N_2(x_0^+) 
    = 0.88525
$$
where $\Cal N_2(x)=\int_{[0,1]^2}|K(x,\tau,\theta)|^2\,d\tau\,d\theta$.
This bound is computed by numerical integration.
To estimate the approximation error, one needs upper bounds of partial
derivatives of $K$. They can be easily obtained using the lower bounds
(\eqMinPsi) and (\eqMinP) of $|\Psi|$ and $|P|$, and upper bounds of 
the derivatives of $\Psi$ obtained from its integral form in (\eqPhi).
For upper bounds of the derivatives of polynomials involved in the definition of $K$
one can use just the sums of upper bounds of monomials.
\qed\enddemo

Replacing the integrals by integral sums, equation
(\eqFredG) can be solved with any given precision.
Then, due to (\eqHg) and (\eqHs) we can numerically compute $H(x)$
using the Cauchy integral
$$
   H(x^3) = \frac{x^2}{2\pi i}\oint_{|t|=1}\frac{g(x,t)\,dt}{t}
      = x^2\int_0^1 \varphi_g(x,\tau)\,d\tau
                                                    \eqno(\eqIntH)
$$
(recall that $\varphi_g(x,\tau):=g(x,e^{2\pi i\tau})$; see (\eqFredTau)).
We can summarize the content of this section as follows
(recall that $f(m,n)$ is the number of primitive lattice triangulations of the
rectangle $m\times n$).

\proclaim{ Proposition \propFred } $\lim_{n\to\infty}f(3,n)^{1/n} = 1/x_0^2$ where:
\roster
\item"$\bullet$"
                 $x_0$ is a unique solution of the equation
                 $H(x^3)=1$ on the interval $[0,x_0^+]$ with $x_0^+=\frac{17}{35}$;
\item"$\bullet$"
                 $H(x)$ is defined via $g(x,t)$ by (\eqIntH)
                 and it is monotone on $[0,x_0^+]$;
\item"$\bullet$"
                 $g(x,t)$ is the solution of the Fredholm equation (\eqFredG)
                 whose ingredients $P$ and $\Psi$ are defined
                 by (\eqDefP) and by $\Psi(x,t)=1-x^2(t-x)\Phi(x,t)$
                 with $\Phi$ defined by (\eqPhi); for any $x\in[0,x_0^+]$
                 the equation (\eqFredG) has a unique solution.
\endroster
\endproclaim

\midinsert
\centerline{\epsfxsize=115mm\epsfbox{prog.eps}}
\botcaption{Figure \figProg}
  Mathematica code for computation of $H(x)$
\endcaption
\endinsert

In Figure~\figProg\ we present a Mathematica function {\,\tt H\,}
which computes $H(x)$ with any given precision.
The approximating error can be estimated using Lemma~\lemFred.
One can check that the functions $P(x,t,u)$ and $\Psi(x,t)$ do 
not vanish when $x<x_0^+$, $|u|=1$, and $\frac{10}{13}<|t|<\frac{13}{10}$.
In Figures~\figPsiImage\ and \figPsiImgII\ we show the image of 
the annulus $\frac{10}{13}<|t|<\frac{13}{10}$ under the mapping $t\mapsto\Psi(x_0,t)$.
Thus we can apply the error estimate (\eqFredErrAA) with $r=10/13$ and hence
$a=-\frac{\log r}{2\pi}=0.04176$.
When estimating the error of $H(x)$ with
$x\approx x_0$, we can set in (\eqFredErrAA)
$$
   C\le 1;\quad
   \tfrac1n\|B\|_1\le 3.05;\quad
   M\le 3910;\quad
   M'\le 94.6;\quad
   M_f\le 258.
$$
Then we obtain the error estimate presented in the last column of Table~\tabErr.
We see that it is reasonably close to the actual error which is given
in the 4th column.

\midinsert
\centerline{\epsfxsize=120mm\epsfbox{psi-image-13.eps}}
\botcaption{Figure \figPsiImage}
    A realistic drawing of the image of the circles $|t|=\tfrac{10}{13}$, $|t|=1$,
    and $|t|=\tfrac{13}{10}$ by the mapping $t\mapsto\Psi(x_0,t)$.
    The left zoom is stretched in the vertical direction.
\endcaption
\bigskip
\bigskip
\centerline{\epsfxsize=100mm\epsfbox{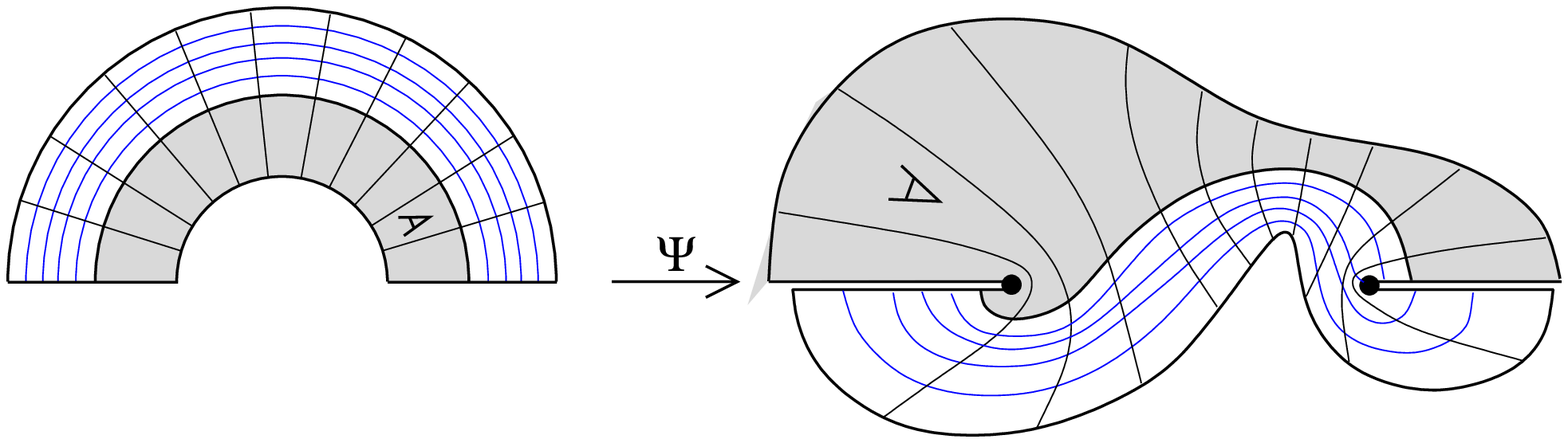}}
\vskip-30pt
\botcaption{Figure \figPsiImgII}
    A schematic drawing of the image of the upper half-annulus
    $\{\frac{10}{13}\le|t|\le \frac{13}{10},\;\Im t\ge0\}$
    by the mapping $t\mapsto\Psi(x_0,t)$. The lower half-annulus
    is mapped symmetrically. The gray part is $|t|<1$.
\endcaption
\endinsert

\midinsert
\hbox to 100mm{\hfill Table \tabErr}
\medskip
\centerline{
\vbox{\offinterlineskip
\def\h {height2pt&\omit&&\omit&&\omit&&\omit&&\omit&\cr}
 \def\t{\times 10} \def\s{\;\;} \def\ss{\s\s} 
\def\q{\quad\;}
\hrule
\halign{&\vrule#&\strut\;\hfil#\hfil\,\cr
\h
&     &&       && time   && $n$-th approx.~of && error            &\cr\h
& $n$ && prec. && (sec.) && $H(x_0)-1$        && estimate         &\cr
\h
 \noalign{\hrule}
\h
& 100  && 24  && 0.299391 && $1.44\t^{-10}$  && $6.95\t^{-4}\ss$ &\cr\h
& 200  && 36  && 6.759046 && $5.01\t^{-22}$  && $5.60\t^{-15}\s$ &\cr\h
& 300  && 48  && 21.77949 && $1.73\t^{-33}$  && $3.39\t^{-26}\s$ &\cr\h
& 400  && 60  && 51.22560 && $6.02\t^{-45}$  && $1.82\t^{-37}\s$ &\cr\h
& 500  && 72  && 115.5499 && $2.09\t^{-56}$  && $9.19\t^{-49}\s$ &\cr\h
& 600  && 84  && 231.5893 && $7.26\t^{-68}$  && $4.45\t^{-60}\s$ &\cr\h
& 700  && 96  && 380.6020 && $2.52\t^{-79}$  && $2.09\t^{-71}\s$ &\cr\h
& 800  && 108 && 608.9937 && $8.78\t^{-91}$  && $9.65\t^{-83}\s$ &\cr\h
& 900  && 120 && 869.7188 && $3.06\t^{-102}$ && $4.38\t^{-94}\s$ &\cr\h
& 1000 && 132 && 1072.923 && $1.06\t^{-113}$ && $1.96\t^{-105}$ &\cr\h
& 1100 && 144 && 1456.021 && $3.72\t^{-125}$ && $8.70\t^{-117}$ &\cr\h
& 1200 && 156 && 1852.763 && $1.29\t^{-136}$ && $3.83\t^{-128}$ &\cr\h
\noalign{\hrule}
}}
}
\endinsert


\bigskip

\head\sectFred. Approximate solutions of Fredholm integral equations with
                analytic kernels
\endhead

\subhead\sectFred.1. Error estimates. Generalities
\endsubhead
The notation in this subsection is independent of the notation in the rest
of the paper.

\proclaim{ Lemma \lemErrL }
Let $f$ be a holomorphic function in a neighborhood of the annulus
$R_1 < |z| < R_2$, and let $f(z)=\sum_{n\in\Z} c_n z^n$ be its Laurent series.
Then, for $R_1<r<R_2$ and for any $n>0$, 
$$
      \Bigg| c_0 - \frac{1}{n}\sum_{k=1}^n f(r\omega^k)\Bigg|
      =\Bigg| \int_0^1 f(re^{2\pi it})\,dt - \frac{1}{n}\sum_{k=1}^n f(r\omega^k)\Bigg|
       \le \frac{M_1 q_1^n}{1-q_1^n} + \frac{M_2 q_2^n}{1-q_2^n}
    \eqno(\eqErrL)
$$ 
where $\omega=e^{2\pi i/n}$, $q_1=R_1/r$, $q_2=r/R_2$,
and $M_j=\max_{|z|=R_j} |f(z)|$ for $j=1,2$.
\endproclaim

\demo{ Proof}
We have
$$
  \sum_{k=1}^n f(r\omega^k)
  = \sum_{k=1}^n\; \sum_{m\in\Z} c_m (r\omega^k)^m
 \quad{and}\quad
   \sum_{k=1}^n \omega^{km}
  = \cases n, &\text{$n$ divides $m$},\\ 0, &\text{otherwise}\endcases
$$
hence the left hand side of (\eqErrL) is equal to
$\big|\sum_{p\in\Z\setminus\{0\}} c_{pn}r^{pn}\big|$ and the coefficients can be estimated
using the Cauchy integrals.
\qed\enddemo

\proclaim{ Lemma \lemTaylor }
Let $h>0$ and $D\subset\R^d$ be a product of segments
$[0,n_1h]\times\dots\times[0,n_dh]$ with positive integers $n_1,\dots,n_d$.
Let $f:D\to\R$ be a function of class $\Cal C^2$
and $M=\max_{i,j}\max_D|\partial_i\partial_j f|$.
Then
$$
    \min_D f \ge \min_{h\Z^d} f - \tfrac18Md^2h^2
$$
where $h\Z^d=\{h\vec n\mid\vec n\in\Z^n\}$.
A similar estimate holds for $\max_D f$.
\endproclaim

\demo{Proof} Induction on $d$. Let the minimum be attained at $x_0\in D$.
If $x_0$ is in the interior of $D$, then we estimate $|f(x)-f(x_0)|$ for the
nearest to $x_0$ grid point $x$ using
the Taylor--Lagrange formula for the second order expansion of
$f(x_0+t(x-x_0))$ at $t=0$.
If $x_0$ is on the boundary of $D$, then we apply the
induction hypothesis to the restriction of $f$ to the facet of $D$
containing $x_0$.
\qed\enddemo


\subhead\sectFred.2. Error estimates for approximate solutions of
                     Fredholm equations
\endsubhead

Let $\varphi:\R\to\C$ be a continuous solution
of the Fredholm integral equation
$$
   \varphi(x) = \int_0^1 K(x,y)\varphi(y)\,dy + f(x)
   \eqno(\eqFredGen)
$$
with analytic complex-valued functions $K$ and $f$ which are (bi)-periodic
  with period $1$, i.e.,
$K(x,y)=K(x+1,y)=K(x,y+1)$ and $f(x)=f(x+1)$.
Assume that $K$ and $f$ extend to complex analytic functions in a
neighborhood of $(D\times\R)\cup(\R\times D_1)$ in $\C^2$ and
in a neighborhood of $D$ in $\C$ respectively where
$$
  D  =\{z\in\C\mid -a  \le \Im z\le a  \}, \quad
  D_1=\{z\in\C\mid -a_1\le \Im z\le a_1\}, \quad
  0<a_1<a.
$$
Let us set
$$
   C=\int_0^1|\varphi(x)|\,dx,
\qquad
   M=\max_{D\times\R}|K|,
\qquad
   M'_1=\max_{\R\times D_1}|K|,
\qquad
   M_f=\max_{D}|f|.
$$

\proclaim{ Lemma \lemFredHolo }
The function $\varphi$ analytically extends to a neighborhood of $D$ and
$$
     M_{\varphi}:=\max_{D_1}|\varphi|
     \le \frac{a(CM+M_f)}{a-a_1}.
     \eqno(\eqFredHolo)
$$
\endproclaim

\demo{ Proof}
For any $(x_0,y_0)\in\R^2$ and any $n$, we have
$$
   \big|\partial_x^n K(x_0,y_0)\big|
   \le\Bigg|\frac{n!}{2\pi i}\int_{|z|=a}\frac{K(z,y_0)\,dz}{z^{n+1}}\Bigg|
   \le \frac{Mn!}{a^n}
$$
and similarly $|f^{(n)}(x_0)|\le M_f n!/a^n$. Then, derivating (\eqFredGen)
$n$ times with respect to $x$, we obtain 
$$
   \big|\varphi^{(n)}(x_0)\big|
    =\Bigg|\int_0^1 \partial_x^n K(x_0,y)\varphi(y)\,dy + f^{(n)}(x_0)\Bigg|
    \le\frac{(CM+M_f)n!}{a^n}.
 \eqno(\eqErrDphi)
$$
Hence the Taylor series of $\varphi$ at $x_0$ converges in the disk $|z-x_0|<a$
and, for $|z-x_0|\le a_1$, we have
$$
   |\varphi(z)|=\Big|\sum_{n\ge0}\frac{\varphi^{(n)}(x_0)}{n!}(z-x_0)^n\Big|
    \le \sum_{n\ge 0}\frac{(CM+M_f)a_1^n}{a^n}
         = \frac{a(CM+M_f)}{a-a_1}
$$
whence the required
bound for $M_\varphi$.
\qed\enddemo

For a positive  integer $n$, let us see what happens
if we replace the integral in (\eqFredGen) by
the $n$-th integral sum. Namely, consider the vectors
$\varphi^{[n]}=(\varphi_1^{[n]},\dots,\varphi_n^{[n]})$,
$f^{[n]}=(f_1^{[n]},\dots,f_n^{[n]})$,
and the $n\times n$ matrix $K^{[n]}=\big(K_{jk}^{[n]}\big)_{jk}$ defined by
$$
  \varphi^{[n]}_j=\varphi(j/n),\quad f_j^{[n]}=f(j/n),\quad
   K_{jk}^{[n]}=\tfrac1n K(j/n,k/n).
$$
Let $\hat\varphi^{[n]}=\big(\hat\varphi^{[n]}_1,\dots,\hat\varphi^{[n]}_n\big)$ be
a solution of the equation
$$
     \hat\varphi^{[n]}=K^{[n]}\hat\varphi^{[n]} + f^{[n]}.
    \eqno(\eqFredAppr)
$$
This equation is a discretization of (\eqFredGen) and it is natural to expect that
$\hat\varphi^{[n]}$ well approximates $\varphi$.
Now, following the approach from [\refKK],
we estimate the rate of the convergence. Our final purpose is to find
a good upper bound for the approximating error
$$
   E_n:=   \Bigg|\int_0^1 \varphi(x)\,dx
                  - \frac1n\sum_{j=1}^n \hat\varphi^{[n]}_j\,\Bigg|.
$$

We define the norms $\|\cdot\|_p$, $1\le p\le\infty$, on $\C^n$ in the usual way.
For a square matrix $A=(a_{jk})_{jk}$ with complex entries we set
$$
  \|A\|_1=\sum_{j,k} |a_{jk}|,
\qquad
    \|A\|_2 
     = \Big(\sum_{j,k}|a_{jk}|^2\Big)^{1/2}.
$$

\proclaim{ Lemma \lemFred }
(a).
Suppose that the matrix $A^{[n]}=I-K^{[n]}$ is invertible and
denote its inverse by $B^{[n]}$.
Then
$$
     E_n
      \le \frac{2a(CM+M_f)r_1^n}{(a-a_1)(1-r_1^n)}\Big(1+\tfrac1n\|B^{[n]}\|_1M'_1\Big),
      \qquad   r_1=e^{-2\pi a_1}.
    \eqno(\eqFredErrA)
$$
If $K$ analytically extends to a neighborhood of $\R\times D$ and $M'=\max_{\R\times D}K$,
$$
     E_n 
      \le 4\pi e(CM+M_f)
     \Big(1+\tfrac1n\|B^{[n]}\|_1 M'\Big)\frac{nar^n}{1-er^n},
      \qquad   r=e^{-2\pi a}.
    \eqno(\eqFredErrAA)
$$
For $n>\alpha_n$, we have
$$
   C\le\frac{\|\hat\varphi^{[n]}\|_1+\alpha_nM_f}{n-\alpha_n}
\qquad\text{where}\qquad
  \alpha_n = \frac{2M'_1a\|B^{[n]}\|_1r_1^n}{(a-a_1)(1-r_1^n)} + \frac{1}{4a}.
    \eqno(\eqFredErrB)
$$

(b).
Suppose that $\|K^{[n]}\|_2=M_2<1$. Then $A^{[n]}$ is invertible and
$\frac1n\|B^{[n]}\|_1\le 1/(1-M_2)$ which implies
in particular that $\alpha_n<\alpha_0$ for some constant
$\alpha_0=\alpha_0(a,a_1,M,M'_1,M_2,M_f)$ and hence
$C$ can be estimated using (\eqFredErrB) for $n>\alpha_0$.

\endproclaim

\demo{ Proof }
(a). Let $J=\int_0^1\varphi(x)\,dx$,
$S=\frac1n\sum_j\varphi(j/n)$,
$\hat S=\frac1n\sum_j\hat\varphi(j/n)$,
$\rho=\varphi^{[n]}-\hat\varphi^{[n]}$,
and $\sigma=A^{[n]}\rho$. In this notation, $E_n=|J-\hat S|$.
We have
$$
 \|\sigma\|_\infty
  = \big\|A^{[n]}\varphi^{[n]}-A^{[n]}\hat\varphi^{[n]}\big\|_\infty
 \overset\text{(\eqFredAppr)}\to=
 \big\|A^{[n]}\varphi^{[n]}-f^{[n]}\big\|_\infty
 =\big\|K^{[n]}\varphi^{[n]}-(\varphi^{[n]}-f^{[n]})\big\|_\infty
$$
By (\eqFredGen), we have
$\varphi^{[n]}_j-f^{[n]}_j=\int_0^1 K(j/n,y)\varphi(y)\,dy$,
and the $j$-th component
of the vector $K^{[n]}\varphi^{[n]}$ is the $n$-th integral sum for this integral.
Hence, applying Lemma~\lemErrL\ to the functions $K(j/n,z(\zeta))\varphi(z(\zeta))$
after the change of variable $\zeta=e^{2\pi iz}$, we obtain
$\|\sigma\|_\infty\le M'_1C_1$ with $C_1=2M_\varphi r_1^n/(1-r_1^n)$ and then
$$
    \|\rho\|_1=\big\|B^{[n]}\sigma\big\|_1
    \le \big\|B^{[n]}\big\|_1\times\|\sigma\|_\infty
    \le M'_1C_1\big\|B^{[n]}\big\|_1.
  \eqno(\eqFredErrRho)
$$
Lemma~\lemErrL\ applied to $\varphi(z(\zeta))$ yields
$|J - S| \le C_1$. We also have $|S-\hat S|\le\frac1n\|\rho\|_1$, hence
$$
  E_n=|J-\hat S|\le|J-S|+|S-\hat S|
    \le C_1 + \tfrac1n\|\rho\|_1
   \le C_1 + \tfrac1n M'_1C_1\|B^{[n]}\|_1
$$
which yields (\eqFredErrA) after applying (\eqFredHolo).
Setting $a_1=a-\frac{1}{2\pi n}$ (hence $r_1=e^{1/n}r$) and $M'_1<M'$
in (\eqFredErrA), we obtain (\eqFredErrAA).

\smallskip
Let us prove (\eqFredErrB).
It is easy to check that
$$
  nC\le\|\varphi^{[n]}\|_1 + \tfrac14\max_{\R}|\varphi'|
    \le\|\hat\varphi^{[n]}\|_1 +\|\rho\|_1  + \tfrac14\max_{\R}|\varphi'|.
$$
Using the estimates (\eqFredErrRho) and (\eqErrDphi) for $\|\rho\|_1$ and $|\varphi'|$
respectively, we obtain
$$
\split
  nC&\le\|\hat\varphi^{[n]}\|_1
  + \frac{2M_1M_\varphi\|B^{[n]}\|_1r_1^n}{1-r_1^n} + \frac{CM+M_f}{4a}
 \overset\text{(\eqFredHolo)}\to\le
  \|\hat\varphi^{[n]}\|_1 + (CM+M_f)\alpha_n.
\endsplit
$$

(b).
Suppose now that $\|K^{[n]}\|_2= M_2<1$.
Then $\|B^{[n]}\|_2=\|(I-K^{[n]})^{-1}\|_2=\|I+K^{[n]}+(K^{[n]})^2+\dots\|_2
\le 1/(1-M_2)$. By Cauchy Inequality we also have
$\|B^{[n]}\|_1\le n\|B^{[n]}\|_2$
\qed\enddemo

\subhead\sectFred.3. A numerical criterion of existence and uniqueness of solutions
\endsubhead
Here we keep the above assumptions about $K(x,y)$ and $f(x)$
except that we no longer assume a priori that equation (\eqFredGen) has
a continuous solution $\varphi$.
Let $\Cal K:\Cal C([0,1])\to\Cal C([0,1])$ be the Fredholm integral operator
with kernel $K(x,y)$, i.e., the operator $\varphi\mapsto\psi$ where
$\psi(x)=\int_0^1 K(x,y)\varphi(y)\,dy$.

\proclaim{ Lemma \lemFredEV } {\rm(cf.~[\refKK, Ch.II, \S1, Eq.~(26)])}.
Suppose that there exists $n$ such that the matrix $I-K^{[n]}$ is invertible
and $\alpha_n<n$ where $\alpha_n$ is
defined in (\eqFredErrB) (note that neither $f$ nor $\varphi$ is used
in the definition of $\alpha_n$). Then $1$ is not an eigenvalue of $\Cal K$
and hence, for any given continuous function $f$, equation (\eqFredGen)
has a unique continuous solution $\varphi$.
\endproclaim

\demo{ Proof } Let $n$ be such that $\alpha_n<n$.
Let us apply Lemma~\lemFred(a) when $f=0$ and hence $\hat\varphi^{[n]}=0$.
Then (\eqFredErrB) reads $C\le 0$ which means that there are no non-zero
solutions of the equation $\Cal K\varphi=\varphi$, i.e. $1$ is not an eigenvalue
of $\Cal K$.
By Fredholm Theorem [\refFred], in this case (\eqFredGen) has
a unique continues solution for any $f$.
\qed\enddemo


\subhead\sectFred.4. Analyticity of solutions with respect to a parameter
\endsubhead 
Let $\Lambda$ be a domain in $\C$ and $U=\{z\in\C\mid -a\le\Im z\le a\}$, $a>0$.
Let $K(\lambda,x,y)$ be an analytic function in a neighborhood of
$\Lambda\times U^2$ in $\C^3$ and $f(\lambda,y)$ be an analytic function
in a neighborhood of $\Lambda\times U$ in $\C^2$.
We assume that $K(\lambda_0,x,y)$
is $(1,1)$-biperiodic and $f(\lambda_0,x)$ is $1$-periodic for any fixed
$\lambda_0\in\Lambda$.

For $\lambda\in\Lambda$,
let $\Cal K_\lambda:\Cal C([0,1])\to\Cal C([0,1])$ be the Fredholm
integral operator $\varphi\mapsto\psi$, $\psi(x)=\int_0^1 K(\lambda,x,y)\varphi(y)dy$.
The next lemma immediately follows from Fredholm's results in
his seminal paper [\refFred]
(a more general fact is proven in [\refTamar]).

\proclaim{ Lemma \lemParam } Suppose that $1$ is not an eigenvalue of $\Cal K_\lambda$
for any $\lambda\in\Lambda$. Then, for any $\lambda\in\Lambda$,
there exists a unique solution $\varphi(\lambda,x)$ of the equation
$$
   \varphi(\lambda,x) = \int_0^1 K(\lambda,x,y)\varphi(\lambda,y)\,dy + f(\lambda,x)
   \eqno(\eqFredParam)
$$
and the function $\varphi(\lambda,x)$ is analytic in a neighborhood of
$\Lambda\times U$.
\endproclaim

\demo{ Proof } By Fredholm's results [\refFred] (see also [\refKhved]),
for any $\lambda\in\Lambda$,
the solution $\varphi(\lambda,t)$ is unique under our assumptions
and it can be written as
$$
     \varphi(\lambda,t)=f(\lambda,t)
     +\int_0^1\frac{D(\lambda,x,y)}{D(\lambda)}f(\lambda,y)\,dy
$$
where
$$
     D(\lambda)=\sum_{n=0}^\infty\frac{(-1)^n A_n(\lambda)}{n!},
\qquad
     D(\lambda,x,y)=\sum_{n=0}^\infty\frac{(-1)^n B_n(\lambda,x,y)}{n!},
    \eqno(\eqFredParamPf)
$$
$$
   A_n(\lambda)=\int_{[0,1]^n} K(\lambda,\bold x,\bold x)\,d\bold x,
   \qquad
    B_n(\lambda,x,y)=\int_{[0,1]^n} K(\lambda,x,\bold x,y,\bold x)\,d\bold x,
$$
$$
   K(\lambda,x_1,\dots,x_n,y_1,\dots,y_n)=\det\big(K(\lambda,x_i,y_j)\big)_{i,j=1}^n.
$$
It is shown in [\refFred] that $D(\lambda)$ does not vanish on $U$
(because $1$ is not an eigenvalue of $\Cal K_\lambda$ for any $\lambda\in\Lambda$).
It is clear that the functions $A_n$ and $B_n$ are analytic in $\Lambda$
and in $\Lambda\times U^2$ respectively and
the Hadamard Inequality $|\det N|\le n^{n/2}\max_{i,j}|N_{ij}|$
implies the upper bounds (cf.~[\refFred, p.~368, line 4]):
$$
     |A_n(\lambda)|\le n^{n/2}M(\lambda)^n,
\qquad
     |B_{n-1}(\lambda,x,y)|\le n^{n/2}M(\lambda)^n
$$
where $M(\lambda)=\sup_{(x,y)\in U^2}|K(\lambda,x,y)|$.
Hence the series (\eqFredParamPf) converge to analytic functions
whence the result.
\qed\enddemo


\head\sectNonprim. Non-primitive lattice triangulations
\endhead

Denote the number of all (not necessarily primitive) lattice
triangulations of the $m\times n$ rectangle by
$f^\np(m,n)$, and set
$$
    c^\np = \lim_{n\to\infty}\frac{\log_2 f^\np(n,n)}{n^2}.
$$

\proclaim{ Proposition \propNP} $c^\np \le 4.735820221...$
\endproclaim

\demo{ Proof } Let $N=n^2$.
Any lattice triangulation can be subdivided up to a primitive
lattice triangulation. Hence a lattice triangulation is completely determined by a
choice of a primitive lattice triangulation and a
set of its edges to be removed.
Let $f_k^\np(n,n)$ be the number
of lattice triangulations of the $n\times n$ square with $k$ interior vertices
and hence with $\approx 3k$ edges. Then
$$
    f_k^\np(n,n) \le \binom{3N}{3k}2^{cN}          \eqno(\eqNPa)
$$
(recall that $2^{cN}$ is a bound for the number of primitive lattice
triangulations).
On the other hand the number of triangulations with vertices in
an arbitrary fixed set of $k$ points on a plane is $O(30^k)$ (see [\refShSh]),
hence
$$
   f_k^\np(n,n) \le \binom{N}{k}30^k.              \eqno(\eqNPb)
$$
Combining (\eqNPa) and (\eqNPb) with Stirling formula, we obtain
$$
   c^\np \le \max_{0\le x\le 1} \min\big(3h(x)+c,h(x)+x\log_2 30\big),  \eqno(\eqNPc)
$$
$$
   h(x)=-x\log_2 x - (1-x)\log_2(1-x).
$$
Using the bound $c\le 4\log_2\frac{1+\sqrt5}2$
(see [\refMVW], [\refWelzl], [\refWelzlTalk]), we obtain the result
(the maximum in (\eqNPc) is attained at $x=0.83206855$).
\qed\enddemo

\Refs
\def\r{\ref}

\r\no\refAnclin
\by E.~Anclin
\paper An upper bound for the number of planar lattice triangulations
\jour  J. Combinatorial Theory, Ser. A \vol 103 \yr 2003 \pages 383--386
\endref

\r\no\refFred
\by  I.~Fredholm
\paper Sur une classe d'\'equations fonctionnelles
\jour Acta Math. \vol 27 \yr 1903 \pages 365--390
\endref

\r\no\refGKZ
\by I.~M.~Gelfand, M.~M.~Kapranov, A.~V.~Zelevinsky
\book Discriminants, Resultants, and Multidimensional Determinants
\publ Birkh\"auser \publaddr Boston \yr 1994
\endref

\r\no\refKZ
\by V.~Kaibel, G.~M.~Ziegler
\paper Counting Lattice Triangulations
\inbook in: C.~D.~Wensley (ed.) Surveys in combinatorics, 2003,
Proc. of the 19th British combinatorial conf., Univ. of Wales, Bangor UK,
June 29 -- July 04
\bookinfo London Math. Soc. Lect. Notes \vol 307 \yr 2003 
\publ Cambridge Univ. Press \publaddr Cambridge
\pages 277--307
\endref

\ref\no\refKK
\by  L.~V.~Kantorovich, V.~I.~Krylov
\book Approximate methods of higher analysis 
\publ Groningen: P. Noordhoff \yr 1958
\endref

\ref\no\refKhved
\by     B.~V.~Khvedelidze
\paper  Fredholm equation
\inbook in: Matematicheskaya Enciklopediya
\ed     I.~M.~Vinogradov 
\publaddr Moscow \yr 1977
\lang Russian
\transl English transl.
Encyclopedia of Mathematics. 
URL: 
http://encyclopediaofmath.org/index.php?title=Fredholm\_equation\&oldid=46977
\endref

\ref\no\refMVW
\by J.~Matou\v sek, P.~Valtr, E.~Welzl
\paper On two encodings of lattice triangulations
\jour manuscript \yr 2006
\endref

\r\no\refOrevkov
\by S.~Yu.~Orevkov
\paper Asymptotic number of triangulations with vertices in $\Bbb Z^2$
\jour  J. Combinatorial Theory, Ser. A \vol 86 \yr 1999 \pages 200--203
\endref

\r\no\refOrKh
\by S.~Yu.~Orevkov, V.~M.~Kharlamov
\paper Asymptotic growth of the number of classes of real plane algebraic
       curves when the degree increases
\jour  Zapiski Nauch. Seminarov POMI \vol 266 \yr 2000 \pages 218--233 
\transl English translation: \jour J. of Math. Sciences \vol 113 \yr 2003 \issue 5 \pages 666--674
\endref

\r\no\refShSh
\by M.~Sharir, A.~Sheffer
\paper Counting triangulations of planar point sets
\jour Electron. J. Combin. \vol 18 \yr 2011 \issue 1 \pages P70:1--74
\endref

\r\no\refTamar
\by J.~D.~Tamarkin
\paper On Fredholm's integral equations, whose kernels are analytic in a parameter
\jour Ann. Math. \vol 28 \yr{1926--1927} \pages 127--152
\endref

\r\no\refWelzl
\by E.~Welzl
\paper The number of triangulations on planar point sets
\inbook In: M.~Kaufmann, D.~Wagner (eds) Graph Drawing. GD 2006
\bookinfo Lecture Notes in Computer Science, vol 4372
\publ Springer \publaddr Berlin, Heidelberg \yr 2007 \pages 1--4
\endref

\r\no\refWelzlTalk
\by E.~Welzl (with J.~Matu\v sek and P.~Valtr)
\paper Lattice triangulations
\jour Talk in Freie Univ. Berlin, November~13, 2006
\endref

\endRefs
\enddocument